\documentclass[webpdf,modern,medium]{oup-authoring-template}
\onecolumn

\usepackage{amsmath,amssymb,mathtools}
\usepackage{mathrsfs}
\usepackage{graphicx}
\usepackage{url}
\usepackage{enumerate}
\usepackage{etoolbox}
\usepackage{tikz}
\usepackage{tikz-cd}
\usetikzlibrary{matrix}

\theoremstyle{thmstyleone}
\newtheorem{theorem}{Theorem}[section]
\newtheorem{lemma}[theorem]{Lemma}
\newtheorem{proposition}[theorem]{Proposition}
\newtheorem{corollary}[theorem]{Corollary}

\theoremstyle{thmstylethree}
\newtheorem{definition}[theorem]{Definition}

\newtheorem{question}[theorem]{Question}

\theoremstyle{thmstyletwo}
\newtheorem{remark}[theorem]{Remark}

\numberwithin{equation}{section}

\newcommand{\Lip}{\operatorname{Lip}}
\newcommand{\Lipz}{\operatorname{Lip}_0}
\newcommand{\LHFz}{\operatorname{HLip}_0}

\newcommand{\BLip}{\operatorname{BLip}}
\newcommand{\BLipz}{\operatorname{BLip}_0}
\newcommand{\HBLipz}{\operatorname{HBLip}_0}
\newcommand{\Blin}{\operatorname{Blin}}
\newcommand{\Fix}{\operatorname{Fix}}
\newcommand{\Coinv}{\operatorname{Coinv}}
\newcommand{\F}{\mathcal{F}}
\newcommand{\Lop}{\mathcal{L}}
\newcommand{\whot}{\widehat{\otimes}_{\pi}}
\newcommand{\bigslant}[2]{{\raisebox{.2em}{$#1$}\left/\raisebox{-.2em}{$#2$}\right.}}

\makeatletter
\patchcmd{\@maketitle}{%
\if@medium%
          \hbox to \textwidth{\raisebox{5pt}[0pt]{%
          \parbox[b]{362pt}{\raggedleft{{\sffamilyfontcnitalic\fontsize{8bp}{12}\selectfont \@journaltitle}, \sffamilyfontcn\@copyrightyear, \thepage--\thelastpage}\\[1pt]
          {\sffamilyfontcn doi: \@DOI}\\[1pt]
          {\ifx\@access\@empty
          \else
          {\sffamilyfontcn \@access}\fi}
          \vskip1pt
          {\ifx\@appnotes\@empty
          \else
          {\sffamilyfontcn \@appnotes}\fi}
            }}%
          \hfill{\color{black!20}\rule{45pt}{55pt}}}
      \else%
}{%
\if@medium%
          \phantom{\hbox to \textwidth{\raisebox{5pt}[0pt]{%
          \parbox[b]{362pt}{\raggedleft{{\sffamilyfontcnitalic\fontsize{8bp}{12}\selectfont \@journaltitle}, \sffamilyfontcn\@copyrightyear, \thepage--\thelastpage}\\[1pt]
          {\sffamilyfontcn doi: \@DOI}\\[1pt]
          {\ifx\@access\@empty
          \else
          {\sffamilyfontcn \@access}\fi}
          \vskip1pt
          {\ifx\@appnotes\@empty
          \else
          {\sffamilyfontcn \@appnotes}\fi}
            }}%
          \hfill{\color{black!20}\rule{45pt}{55pt}}}}%
      \else%
}{}{\PackageWarning{manuscript}{Could not hide the Modern Medium opening-page banner}}

\patchcmd{\@maketitle}{%
    \if@medium\noindent\rule{\textwidth}{1\p@}\else%
}{%
    \if@medium\noindent\phantom{\rule{\textwidth}{1\p@}}\else%
}{}{\PackageWarning{manuscript}{Could not hide the Modern Medium opening-page rule}}
\makeatother

\makeatletter
\def\ps@opening{%
  \let\@oddhead\@empty
  \let\@evenhead\@empty
  \let\@oddfoot\@empty
  \let\@evenfoot\@empty
}
\makeatother

\begin{document}

\makeatletter
\renewenvironment{proof}[1][\proofname]{\par\removelastskip%
  \pushQED{\ensuremath{\qed}}%
  \normalfont \topsep7.5\p@\@plus7.5\p@\relax%
  \trivlist%
  \item[\hskip\labelsep\itshape #1.\ ]\ignorespaces%
}{%
  \popQED\endtrivlist\@endpefalse%
}
\makeatother

\journaltitle{}
\DOI{}
\copyrightyear{}
\pubyear{}
\vol{}
\issue{}
\access{}
\appnotes{Paper}
\firstpage{1}

\title[Tensor preduals in two-Lipschitz spaces]
{Tensor preduals, positive homogeneity and harmonicity in two-Lipschitz spaces}

\author[1,$\ast$]{Arindam Mandal}

\address[1]{%
\orgdiv{School of Mathematical Sciences},
\orgname{National Institute of Science Education and Research (NISER),
Homi Bhabha National Institute},
\orgaddress{%
Jatni,
Odisha,India}}

\corresp[$\ast$]{Corresponding author,
\href{mailto:munnainmath@gmail.com}{Email: munnainmath@gmail.com}}

\abstract{We study a two-variable framework for distinguished subspaces of
two-Lipschitz mappings through group actions and tensor preduals. We show that the
ambient spaces are generally large: when $M$ and $N$ are infinite metric spaces,
$\F(M)\widehat{\otimes}_{\pi}\F(N)$ contains an isomorphic copy of
$\ell_1(\mathbb N\times\mathbb N)$, and consequently
$\BLip_0(M,N)$ contains an isomorphic copy of
$\ell_\infty(\mathbb N\times\mathbb N)$. This motivates the study of
smaller subspaces selected by symmetries. We consider group actions by
similarities and the associated character-equivariant two-Lipschitz
mappings. These classes arise as fixed-point spaces of induced
representations, are $weak^*$ closed, and admit quotient tensor
preduals determined by the corresponding coinvariant relations.\\
Two applications illustrate the framework. For the product dilation
action on Banach spaces, character-equivariance is exactly positive
homogeneity in each variable. The resulting positively homogeneous
two-Lipschitz space admits a canonical tensor-predual description and
two different projection procedures. Radial extension yields an
explicit bounded projection for arbitrary Banach-valued mappings,
whereas amenable averaging gives a contractive, and in the nontrivial
case norm-one, projection for dual-valued mappings. Concrete examples
show that the dimension of this subspace can vary substantially with
the underlying spaces.\\
The second application concerns harmonic two-Lipschitz mappings on
finitely generated groups. Using mixed discrete derivatives and
commuting averaging operators, we characterize harmonicity in each
variable and relate it to a joint eigenspace description. A normalized
translation action yields an invariant $weak^*$-closed harmonic
subspace whose predual factors as the projective tensor product of the
corresponding one-variable harmonic quotient preduals. This
factorization gives finite-dimensional consequences and leads to
questions concerning complementability of the harmonic relation
spaces and the relation between separate harmonicity and harmonicity
with respect to the product measure.}

\keywords{Lipschitz map; Lipschitz-free space; two-Lipschitz map; group action by similarities; coinvariant quotient; positive homogeneity; harmonic two-Lipschitz map; joint eigenspace.}

\maketitle
\begingroup
\renewcommand{\thefootnote}{}
\footnotetext{%
\textit{2020 Mathematics Subject Classification.} Primary 46B20; Secondary 46M05, 43A07.\\}

\endgroup

\section{Introduction}
The main motivation of this article is to study the interplay between group actions and \emph{two-Lipschitz mappings}. The corresponding character-equivariant classes form distinguished \(weak^*\)-closed subspaces of the associated two-Lipschitz spaces. These subspaces are often easier to study than the much larger ambient spaces. In some elementary cases, their dimensions can also be determined explicitly.

Given a pointed metric space $(M,d_M,0)$ with distinguished point $0$,
there is a Banach space $\F(M)$, unique up to linear isometry, known as
the \emph{Lipschitz-free space} over $M$. The construction is closely
related to the classical Arens--Eells space introduced in
\cite{OEUATS}; see also \cite{LFBS,LA} for modern accounts of
Lipschitz-free spaces.  A map $f: M \to E$ from the metric space $M$ to a real Banach space $E$ is said to be Lipschitz if $\Lip(f)
=
\sup_{ x\neq y}
\frac{\|f(x)-f(y)\|}{d(x,y)} < \infty
$.
We denote by $\operatorname{Lip}_0(M, E)$ the Banach space of Lipschitz functions on $M$ vanishing at $0$, equipped
with the norm $\Lip(\cdot)$. Then the \emph{Lipschitz-free space} over $M$ is formally defined as
$$\mathcal{F}(M)
=
\overline{\operatorname{span}}
\{\delta_M(m):m\in M\}
\subseteq \operatorname{Lip}_0(M,\mathbb{R})^*,$$
where the closure is taken in the operator norm and for each $m\in M$, $\delta_M(m)\in\operatorname{Lip}_0(M,\mathbb{R})^*$
be the evaluation functional given by
$ \delta_M(m)\left(f\right)=f(m)$ for all
$f\in\operatorname{Lip}_0(M,\mathbb{R})$. The map $\delta_M: M \to \F(M)$ preserving the point $0$ is a nonlinear isometric embedding.
 
The \emph{Lipschitz-free space} provides a natural bridge between the nonlinear metric geometry of a pointed metric space and the linear geometry of an associated Banach space. This connection is expressed through the following universal property:

For any $f \in \Lip_0(M,E)$,
there exists a unique bounded linear operator
$\widehat{f}\colon\mathcal{F}(M)\to E$ such that
$\widehat{f}\circ\delta_M=f$.
Moreover,
$\|\widehat{f}\|=\operatorname{Lip}(f)$.
Thus, the following diagram commutes:
\[
\begin{tikzcd}[column sep=large,row sep=large]
M \arrow[r,"f"] \arrow[dr,"\delta_M"']
  & E \\
  & \mathcal{F}(M) \arrow[u,"\widehat{f}"']
\end{tikzcd}
\]
Applying the universal property to real-valued
maps yields the space $\mathcal{F}(M)$ is a canonical predual of
$\operatorname{Lip}_0(M, \mathbb{R})$. We refer \cite{LA,LFBS} for the basics and preliminaries of the Lipschitz space theory.

\medskip
\noindent\textbf{Notation and conventions.}
Throughout, all metric spaces are assumed to be pointed, and their distinguished point will always be denoted by \(0\), irrespective of the particular metric space under consideration. For a normed linear space, the distinguished point is understood to be its additive identity. All normed linear spaces are over the real field. Unless otherwise specified, \(G\) and \(H\) denote groups with the identity elements $e_G,e_H$, respectively; when finite generation is required, this will be stated explicitly. The notations
$\mathcal{L}(X,Y),~
X\cong Y
~\text{and}~
X\widehat{\otimes}_{\pi}Y
,~\Blin(X,Y;E)$
denote, respectively, the space of bounded linear maps between the Banach spaces \(X\) and \(Y\); $X,Y$ are linearly isometrically isomorphic, and the completion of their projective tensor product, and the space of bounded bilinear maps between $X \times Y$ and Banach space $E$. By \(\mathbb{S}_X\) we 
denote the unit sphere of a normed linear space \(X\). We use \(\mathbb{}\) to denote the real field. For any set \(A\), \(\mathbb{R}^A\) denotes the space of all maps from \(A\) into \(\mathbb{R}\). By $\operatorname{Isom}(M)$ we consider the group of isometries on metric space $M$ and $GL(X)$ signifies the group of invertible bounded linear maps on the Banach space $X$.

The two-variable counterpart of the Lipschitz theory is provided by
two-Lipschitz mappings. For the basics of the two-variable theory, we refer \cite{TLOI}. For pointed metric spaces $M$ and $N$, and normed linear space $E$,  the space
$\BLip_0(M,N;E)$ consists of those mappings $T:M\times N\to E$ vanishing
on the coordinate axes for which the mixed difference
$T(x,y)-T(x',y)-T(x,y')+T(x',y')$
is controlled by a constant multiple of
$d_M(x,x')d_N(y,y')$. The smallest such constant gives rise to a norm, denoted by $\BLip(\cdot)$ on $\BLipz(M,N;E)$. As recalled in Section~2, every
$T\in\BLip_0(M,N;E)$ admits a canonical (bi-)linearization, and hence
\[
\BLip_0(M,N;E)
\cong
\Blin(\F(M),\F(N);E)
\cong
\mathcal L\bigl(\F(M)\widehat\otimes_\pi\F(N),E\bigr)~~\mbox{(see \cite{AM1})}.
\]
  In particular,
$\BLip_0(M,N):=\BLip_0(M,N;\mathbb{R})
\cong
\bigl(\F(M)\widehat\otimes_\pi\F(N)\bigr)^*$.
Thus, in the two-variable setting, the projective tensor product
$\F(M)\widehat\otimes_\pi\F(N)$ plays the role occupied by a single
\emph{Lipschitz-free space} in one-variable Lipschitz theory.

This immediately suggests studying distinguished $weak^*$-closed
subspaces of $\BLip_0(M,N)$ through suitable quotients of its tensor
predual. The motivation is reinforced by the size of the ambient
spaces. For an infinite pointed metric space $M$, the space
$\Lip_0(M,\mathbb R)$ contains an isomorphic copy of $\ell_\infty$, equivalently,
$\F(M)$ contains a complemented copy of $\ell_1$; see
\cite{OTSOLFS}. This has motivated the study of smaller quotient
preduals associated with $weak^*$-closed subspaces of Lipschitz spaces;
see, for instance, \cite{LHF,AM2,DouchaIPC}. The same phenomenon persists
in two variables: if $M$ and $N$ are infinite, then
$\F(M)\widehat\otimes_\pi\F(N)$
contains an isomorphic copy of
$\ell_1(\mathbb N\times\mathbb N)$, and consequently
$\BLip_0(M,N)$ contains an isomorphic copy of
$\ell_\infty(\mathbb N\times\mathbb N)$; see
Proposition~\ref{prop:linfty-embedding}. Hence the natural tensor
predual is again, in general, a large Banach space.

The purpose of the present paper is to develop a systematic mechanism
for producing and describing smaller, structurally distinguished
subspaces of $\BLip_0(M,N)$. 
A convenient class of actions for the present purpose is given by
\emph{actions by similarities}. If a group $G$ acts on a nontrivial metric
space $M$ by similarities, then each group element scales all distances by
a positive factor. These scaling factors are multiplicative under the group
operation and therefore determine a positive character
$$\chi_G:G\longrightarrow (0,\infty)~\mbox{
such that}~
d_M(gx,gx')=\chi_G(g)\,d_M(x,x'),~ \forall ~g\in G,\ x,x'\in M.$$
Thus actions by similarities may be viewed as a natural enlargement of
isometric actions, the latter corresponding to the special case
$\chi_G\equiv 1$. We recall the precise formulation and its elementary
proof in Section~\ref{subsec:act by sim}; see also the standard discussion
of group actions and equivariance in \cite{DAS, tomDieck}, \cite[p.~694]{MSLOGFHM}.

Now suppose that $G$ and $H$ act by similarities on pointed metric spaces
$M$ and $N$, with associated similarity characters $\chi_G$ and
$\chi_H$, respectively. Their product action on product metric space $M\times N$ carries the
product character
$\chi(g,h)=\chi_G(g)\chi_H(h)$.
Since the defining two-Lipschitz estimate is multiplicative in the two
metric variables, this product character is exactly the natural scaling
factor for the induced action on $\BLip_0(M,N)$. This leads to the class
of \emph{character-equivariant two-Lipschitz mappings} (see \ref{subsec:character-equivariant-two-Lipschitz}) and, more generally, to
the invariant-subspace framework developed in Section \ref{sec:invariant-two-lipschitz}.
The product action above induces compatible isometric representations on
both $\BLip_0(M,N)$ and its tensor predual
$\F(M)\widehat\otimes_\pi\F(N)$.
The duality between these representations allows invariant
pointwise-closed subspaces to be described through annihilators, closed
relation spaces and quotient preduals. In particular,
the space of \emph{character-equivariant two-Lipschitz mappings} appear as fixed-point spaces of the representation,
while the corresponding preduals are naturally realized as coinvariant
tensor quotients.


\medskip
\noindent\textbf{Main results.}
The first main result of the paper gives such a description at the
level of invariant pointwise-closed subspaces.  Let
$\mathcal S\subseteq\BLip_0(M,N)$ be a linear subspace which is closed
for pointwise convergence and invariant under the induced product
action.  We prove that $\mathcal S$ is automatically $weak^*$ closed.
If
$ K=\mathcal S_\perp
   \subseteq
   \F(M)\widehat\otimes_\pi\F(N)$,
then $K$ is invariant under the corresponding predual action and
$\mathcal S
   \cong
   \left(
     \bigslant{\F(M)\widehat\otimes_\pi\F(N)}{K}
   \right)^*$.  Moreover, the relation space $K$ is generated, after
taking norm closure, by finite tensor relations.  In this sense,
pointwise linear constraints on two-Lipschitz mappings correspond
precisely to norm-closed relations in the tensor predual.

A particularly important case occurs when $\mathcal S$ is a fixed-point
space.  If $U:=(U_{g,h})_{g \in G, h \in H}$ denotes the induced representation on
$\BLip_0(M,N)$ and $W:=(W_{g,h})_{g \in G,h \in H}$ the corresponding representation on
$\F(M)\widehat\otimes_\pi\F(N)$, then we identify
$\left\{T \in \BLipz(M,N): U_{g,h}T=T\right\}=:\Fix(U)
   =
   \Coinv(W)^\perp$,
where $\Coinv(W)$ is the norm-closed span of the differences
$W_{g,h}u-u$.  Consequently,
$$\Fix(U)
   \cong
   \left(
      \bigslant{\F(M)\widehat\otimes_\pi\F(N)}
           {\Coinv(W)}
   \right)^*.$$
Thus, fixed points in the two-Lipschitz space are dual to coinvariants
of the tensor predual.  For actions by similarities, the fixed-point
space is exactly the \emph{space of character-equivariant two-Lipschitz maps}.
When $G$ and $H$ are amenable, invariant means further yield
contractive projections onto these fixed-point spaces. The averaging
argument also extends naturally to dual-valued two-Lipschitz mappings.
Related norm-one projection phenomena for one-variable Lipschitz
mappings with dual-valued codomains were considered in \cite{AM3}.
For the basics of invariant means, we refer to \cite{GNFA}.

We give two rather different applications of this general framework.

The first concerns positive homogeneity.  Let $X$ and $Y$ be Banach
spaces and let $(0,\infty)^2$ act on $X\times Y$ by independent
dilations.  The corresponding character is
$ \chi(a,b)=ab$,
and the associated character-equivariant maps are precisely the
two-Lipschitz mappings which are positively homogeneous in each
variable.  We show that their canonical predual is obtained by
factoring out the relations
$ \delta_{ax}-a\delta_x~\text{and}~
   \delta_{by}-b\delta_y$ for $a,b \in (0,\infty), x \in X, y \in Y$
in the two free-space factors.  If
$\F^{ph}(X)$ and $\F^{ph}(Y)$ denote the resulting one-variable
quotient preduals, then
$\BLip^{ph}_0(X,Y)
   \cong
   \bigl(
      \F^{ph}(X)\widehat\otimes_\pi\F^{ph}(Y)
   \bigr)^*$.
  The same quotient also gives the vector-valued
linearization of positively homogeneous two-Lipschitz mappings.

This quotient description is complemented by a geometric one.
Restriction of two-Lipschitz map to the product of the unit spheres determines a
positively homogeneous two-Lipschitz map completely, and radial
extension gives the inverse correspondence.  This produces an explicit
bounded projection from $\BLip_0(X,Y;E)$ onto
$\BLip^{ph}_0(X,Y;E)$ for an arbitrary Banach space $E$ (see Proposition \ref{prop:radial-projection}).  On the other
hand, since the dilation group is amenable, invariant averaging yields
a norm-one projection whenever the range space is a dual space (see Theorem \ref{thm:amenable-projection}).  The
two constructions are of a different nature: the radial projection is
explicit and geometric, independent of any condition on the co-domain space, whereas the invariant-mean projection is
contractive but generally nonconstructive, dependent on the duality of the co-domain space.  We also obtain concrete
dimension results.  In particular,
$\BLip^{ph}_0(\mathbb R,\mathbb R;\mathbb R)$ is four-dimensional,
while
$\BLip^{ph}_0(\mathbb R^2,\mathbb R^2;\mathbb R)$ is
infinite-dimensional.

The second application concerns harmonic mappings on finitely generated
groups.  This part is motivated by the discrete differential
description of ordinary Lipschitz functions on groups (see \cite[Section~1.7.2, p.~31]{Yadin}).  For a
finitely generated group, the Lipschitz seminorm may be detected by
right derivatives along a finite symmetric generating set.  We develop
the corresponding two-variable construction by introducing mixed right
derivatives and a discrete two-gradient $\nabla^{(2)} :\left\{T: G\times H \to \mathbb R\right\} \to \left(\mathbb R^{S_G \times S_H}\right)^{G\times H}$, where $G,H$ be two finitely generated groups with the symmetric generating sets $S_G,S_H$, respectively.  For 
$T\in \BLipz(G,H)$ we prove the exact identity
$\BLip(T)=\|\nabla^{(2)}T\|_\infty$.
The class of maps with bounded two-gradient is independent of the
particular finite symmetric generating sets, up to equivalence of the
corresponding norms. We also characterize the functions with vanishing mixed two-gradient:
$\nabla^{(2)}T=0$ precisely when
$T(x,y)
=
T(x,e_H)+T(e_G,y)-T(e_G,e_H)$,
or, equivalently, when $T(x,y)=f(x)+g(y)$ for suitable functions
$f:G\to\mathbb R$ and $g:H\to\mathbb R$. This discrete differential
description also provides the starting point for the harmonic theory
developed in the following subsection.

Further, let $\mu$ and $\nu$ be finitely supported probability measures
on $G$ and $H$, respectively. The discrete differential viewpoint above
provides a natural motivation for harmonicity. Namely, after introducing
the right derivatives in the two variables, one may average them with
respect to $\mu$ and $\nu$. Requiring these weighted averages to vanish
leads to the corresponding notions of $\mu$-harmonicity in the first
variable and $\nu$-harmonicity in the second. We therefore study
two-Lipschitz mappings which are harmonic separately in the two
variables; see Section~\ref{harmonicity}.

This harmonic construction is also connected with the preceding
equivariant theory through a spectral viewpoint. For suitable group
actions, \emph{character-equivariant two-Lipschitz mappings} are simultaneous
eigenvectors of the corresponding averaging operators, whereas harmonic
two-Lipschitz mappings are precisely those fixed by both operators. The
two averaging operators commute, and separate harmonicity implies
harmonicity with respect to the product measure $\mu\otimes\nu$.

Since ordinary translations do not preserve the normalization on the
coordinate axes, we introduce a recentered translation action on
$\BLip_0(G,H)$. This action is isometric and preserves the two-gradient.
The space
$\HBLipz^{\mu,\nu}(G,H)$
of two-Lipschitz mappings which are harmonic in each variable is
invariant and pointwise closed, and hence weak$^*$ closed.

The harmonic subspace also admits a tensor-predual description.
Writing
$\F_\mu^h(G)=\F(G)/N_\mu(G),
   ~
   \F_\nu^h(H)=\F(H)/N_\nu(H)$,
for the one-variable harmonic quotient preduals, we obtain
$\HBLipz^{\mu,\nu}(G,H)
   \cong
   \bigl(
      \F_\mu^h(G)
      \widehat\otimes_\pi
      \F_\nu^h(H)
   \bigr)^*$ (see Subsection \ref{har ten prdu}).  This immediately transfers finite-dimensional
information from the one-variable harmonic theory to the two-variable
setting.  In particular, whenever
$\LHFz^{\mu}(G)$ and $\LHFz^{\nu}(H)$ are finite-dimensional,
$\dim\HBLipz^{\mu,\nu}(G,H)
   =
   \dim\LHFz^{\mu}(G)\,
   \dim\LHFz^{\nu}(H)$.

A one-variable character-equivariant framework was considered in
\cite{AMA}. In the present two-variable setting, however, the natural
predual is the projective tensor product
$\F(M)\widehat{\otimes}_{\pi}\F(N)$.
Accordingly, relations in the two variables interact through the tensor
structure and give rise to product coinvariants and tensor factorizations
of quotient preduals. This provides the structural basis for the
positively homogeneous and harmonic theories developed below; the latter,
in particular, is not considered in \cite{AMA}.

\section{Preliminaries}
\subsection{Group-theoretic notions} \label{subsec: group theo notion}

We begin by recalling some standard terminology concerning group actions;
see \cite[Definition~2.2.1]{DAS}, \cite{tomDieck}, \cite{EK}. Let $G$ be a group with identity
element $e_G$, and let $S$ be a nonempty set. An action of $G$ on $S$ is
determined by a group homomorphism
$\rho\colon G\longrightarrow \operatorname{Sym}(S)$,
where $\operatorname{Sym}(S)$ denotes the group of all bijections of
$S$. As usual, we write
$gx:=\rho(g)(x)$, for all $g\in G,\ x\in S$.
A set equipped with such an action is called a $G$-set. If $M$ and $N$
are $G$-sets, a map $f\colon M\to N$ is said to be $G$-equivariant if for all $g\in G,\ x\in M, ~ ~f(gx)=gf(x)$.

In particular, suppose that $M$ and $N$ are metric spaces and that $G$
acts on them through homomorphisms
\[
\alpha\colon G\longrightarrow \operatorname{Isom}(M)
\quad\text{and}\quad
\beta\colon G\longrightarrow \operatorname{Isom}(N).
\]
Then a map $f\colon M\to N$ is $(\alpha,\beta)$-equivariant precisely
when
$f(\alpha(g)x)=\beta(g)f(x)$ for all $ g\in G,\ x\in M$.
This point of view is standard in the study of equivariant Lipschitz
mappings; see, for example, \cite[p.~694]{MSLOGFHM}.

For our purposes, we shall focus on a particular class of group actions,
namely actions by similarities. The preceding formulation of
$(\alpha,\beta)$-equivariance then provides a natural motivation for the
notion of a \emph{character-equivariant two-Lipschitz mapping}, introduced
in Definition~\ref{def:character-equivariant-two-Lipschitz}. As will be seen throughout the paper,
this notion provides a common framework for several distinguished
classes of two-Lipschitz mappings.

\subsection{Group actions by similarities} \label{subsec:act by sim}

Let $(M,d)$ be a nontrivial metric space. We denote by
$\operatorname{Sim}(M)$ the group of all similarities of $M$, that is,
\[
\operatorname{Sim}(M)
:=
\left\{
T\colon M\to M \text{ bijective} :
d(Tx,Ty)=c_Td(x,y)
\text{ for some }c_T>0
\right\}.
\]
Since $M$ is nontrivial, the constant $c_T$ is uniquely determined by
$T$ and will be referred to as the \emph{similarity ratio} of $T$. The map
$\lambda\colon \operatorname{Sim}(M)\longrightarrow (0,\infty),~
\lambda(T)=c_T$,
is a group homomorphism. Indeed, for $S,T\in\operatorname{Sim}(M)$ and
$x,y\in M, ~ ~
d(STx,STy)
=
c_S d(Tx,Ty)
=
c_Sc_Td(x,y),$
and hence
\[
\lambda(ST)=\lambda(S)\lambda(T).
\]

Suppose now that $G$ acts on $M$ by similarities; equivalently, the
action is given by a homomorphism
$
\rho\colon G\longrightarrow \operatorname{Sim}(M)$.
The composition
$\chi:=\lambda\circ\rho\colon G\longrightarrow(0,\infty)$
is then a positive character of $G$, and $d(gx,gy)=\chi(g)d(x,y),$ for all $ g\in G,\ x,y\in M$.
We shall refer to $\chi$ as the \emph{similarity character} associated
with the action.

The usual isometric actions occur as the special case
$\chi\equiv 1$. Indeed, in this case $d(gx,gy)=d(x,y),$ for all $ g\in G,\ x,y\in M$,
so that the action takes values in $\operatorname{Isom}(M)$.

The preceding discussion may be summarized as follows.

\begin{proposition}\label{prop:similarity-character}
Let $(M,d)$ be a nontrivial metric space and let $G$ act on $M$. Then
the following assertions are equivalent:
\begin{itemize}
\item[(i)] the action of $G$ on $M$ is by similarities;
\item[(ii)] there exists a positive character
$\chi\colon G\longrightarrow(0,\infty)$
such that
\[
d(gx,gy)=\chi(g)d(x,y), ~\mbox{for any}~  g\in G,\ x,y\in M.
\]
\end{itemize}
Moreover, whenever these conditions hold, the character $\chi$ is
uniquely determined by the action.
\end{proposition}

\begin{proof}
Suppose first that the action is by similarities, and let
$\rho:G\to\operatorname{Sim}(M)$ denote the corresponding homomorphism. Then $\chi=\lambda\circ\rho$ is a positive character and,
by the definition of the similarity ratio, $d(gx,gy)=\chi(g)d(x,y)$ for every $g\in G$ and $x,y\in M$.

Conversely, suppose that such a positive character $\chi$ exists. Then for
each $g\in G$, the map $x\mapsto gx$ is a bijection and satisfies
$d(gx,gy)=\chi(g)d(x,y),$ for all $x,y\in M$.
Thus $g$ acts as a similarity with similarity ratio $\chi(g)$.
Uniqueness of $\chi$ follows from the assumption that $M$ contains at
least two distinct points.
\end{proof}

\noindent\textbf{Convention.}
Throughout the remainder of the paper, unless explicitly stated otherwise,
whenever a positive character $\chi\colon G\to(0,\infty)$ is considered
in connection with an action of $G$ on a metric space $(M,d)$, it is
understood that the action is by similarities and that $\chi$ is the
associated similarity character; that is,
$d(gx,gy)=\chi(g)d(x,y)$ for all $g\in G,~x,y\in M$. For any two pointed metric spaces $M,N$, we consider the product metric on $M \times N$ with the distinguished point $(0_M,0_N)$.

\subsection{Character-equivariant two-Lipschitz mappings}
\label{subsec:character-equivariant-two-Lipschitz}

Let $(M,d_M,0)$ and $(N,d_N,0)$ be nontrivial pointed metric spaces.
Suppose that a group $G$ acts on $M$ and a group $H$ acts on $N$ by
similarities, fixing the distinguished points. Let
$\chi_G:G\longrightarrow(0,\infty)~\text{and}~\chi_H:H\longrightarrow(0,\infty)$
denote the corresponding similarity characters. Thus
$d_M(gx,gx')
=
\chi_G(g)d_M(x,x')$
for all $g\in G$ and $x,x'\in M$, and
$d_N(hy,hy')
=
\chi_H(h)d_N(y,y')$
for all $h\in H$ and $y,y'\in N$.

 A map $T:M\times N\to E$ which vanishes
on the coordinate axes and satisfies
$$\|T(x,y)-T(x',y)-T(x,y')+T(x',y')\|\leq C\,d_M(x,x')d_N(y,y')
$$
for some $C\geq0$, is said to be \emph{two-Lipschitz map}.  The least such constant is denoted by $\BLip(T)$.
We use the canonical isometric identification
\begin{equation}\label{eq:BLip-predual}
\BLip_0(M,N)
:=\BLip_0(M,N;\mathbb R)\cong
\bigl(\F(M)\whot\F(N)\bigr)^*,
\end{equation}
where
$T_L(\delta_x\otimes\delta_y)=T(x,y)$.
For convenience, we put
$\mathcal Z(M,N)=\F(M)\whot\F(N)$
and
$\mathscr M_2(M,N)
=
\operatorname{span}\{\delta_x\otimes\delta_y:x\in M,\ y\in N\}$.
Then $\mathscr M_2(M,N)$ is norm dense in $\mathcal Z(M,N)$.

\medskip
\noindent\textbf{The product character and representation.}
The product action of $G\times H$ on $M\times N$ is given by
$(g,h)\cdot(x,y)=(gx,hy)$.
Associated with this action we consider the product character
$$\chi:G\times H\longrightarrow(0,\infty),
~\mbox{given by}~\chi(g,h)=\chi_G(g)\chi_H(h).$$
Since $\chi_G$ and $\chi_H$ are group homomorphisms, so is $\chi$. This is the usual product construction for characters; see
\cite[Chapter~3, \S3.2]{SerreRep}.

Let $E$ be a Banach space. The character $\chi$ gives rise to the
natural representation
$\rho_\chi:G\times H\longrightarrow GL(E)$, given by $\rho_\chi(g,h)=\chi(g,h)I_E$.
It is therefore natural to consider two-Lipschitz mappings which are
equivariant with respect to the product similarity action of
$G\times H$ on $M\times N$ and the representation $\rho_\chi$ on $E$.
The corresponding equivariance condition is
$T(gx,hy)
=
\rho_\chi(g,h)T(x,y)
=
\chi(g,h)T(x,y)$.
This motivates the following definition.

\begin{definition}\label{def:character-equivariant-two-Lipschitz}
A map
$T\in\BLip_0(M,N;E)$
is said to be \emph{$\chi$-equivariant} if
$T(gx,hy)
=
\chi(g,h)T(x,y)$
for all $g\in G$, $h\in H$, $x\in M$, and $y\in N$.

We denote the space of all such \emph{character equivariant two-Lipschitz mappings} by
$\BLip_0^{G\times H,\chi}(M,N;E)$.
In the scalar-valued case, we write simply
$\BLip_0^{G\times H,\chi}(M,N)
:=
\BLip_0^{G\times H,\chi}(M,N;\mathbb R)$.
\end{definition}

Thus $\chi$-equivariance is precisely the usual equivariance condition
associated with the product similarity action of $G\times H$ on
$M\times N$ and the scalar representation
$(g,h)\longmapsto \chi(g,h)I_E$
on $E$.

In particular, if both actions are isometric, then
$\chi_G\equiv1 ~\text{and}~
\chi_H\equiv1$,
and hence
$\chi\equiv1$.
In this case, the $\chi$-equivariance condition reduces to
$T(gx,hy)=T(x,y)$
for all $g\in G$, $h\in H$, $x\in M$, and $y\in N$. Therefore, $T$ is
constant on the $(G\times H)$-orbits of $M\times N$. Thus the case
$\chi\equiv1$ may be regarded as the unweighted special case of the
character-equivariant framework.

\medskip
\noindent\textbf{The dilation action.}
The product dilation action provides one of the examples of the
preceding construction. Let $X$ and $Y$ be Banach spaces and we consider the groups
$G=H=(0,\infty)$,
with multiplication, acting on $X$ and $Y$, respectively, by
$a\cdot x=ax, ~\mbox{and}~
b\cdot y=by$.
The corresponding similarity characters are $\chi_G(a)=a, ~\chi_H(b)=b$,
and hence the associated product character is $\chi(a,b)=ab$.
Thus, for $T\in\BLip_0(X,Y;E)$, the $\chi$-equivariance condition
takes the form
$T(ax,by)=ab\,T(x,y),
~\mbox{for any}~ a,b>0$.

This naturally leads to the following class of two-Lipschitz mappings.

\begin{definition}\label{def:positively-homogeneous-two-Lipschitz}
Let $X$, $Y$, and $E$ be Banach spaces. A map $T\in\BLip_0(X,Y;E)$
is said to be \emph{positively homogeneous in each variable} if
$T(ax,y)=aT(x,y)~~\text{and}~~
T(x,by)=bT(x,y)$
for all $a,b\geq0$, $x\in X$, and $y\in Y$.

We denote the space of all such \emph{positively homogeneous two-Lipschitz mappings} by
$\BLip_0^{ph}(X,Y;E)$.
In the scalar-valued case, we write
$\BLip_0^{ph}(X,Y)
:=
\BLip_0^{ph}(X,Y;\mathbb R)$.
\end{definition}

For the product dilation action, $\chi$-equivariance and positive
homogeneity are equivalent. Indeed, if
$T(ax,by)=ab\,T(x,y)$
for all $a,b>0$, then taking $b=1$ and $a=1$, respectively, gives
$T(ax,y)=aT(x,y)~\text{and}~T(x,by)=bT(x,y)$.
The case of the scalar $0$ follows from the defining conditions $T(0,y)=T(x,0)=0$.
Conversely, positive homogeneity in each variable immediately yields
$T(ax,by)=ab\,T(x,y)$.
Consequently,
$\BLip_0^{(0,\infty)^2,\chi}(X,Y;E)
=
\BLip_0^{ph}(X,Y;E)$,
where
$\chi(a,b)=ab$.

\section{Invariant pointwise-closed subspaces of two-Lipschitz spaces}
\label{sec:invariant-two-lipschitz}

Motivated by the viewpoint developed for invariant pointwise-closed
subspaces of Lipschitz spaces in \cite{DouchaIPC}, we develop
here a two-variable tensor analogue for \emph{two-Lipschitz mappings}.  The
natural predual in this setting is the projective tensor product of two
\emph{Lipschitz-free spaces}.  We first describe invariant pointwise-closed
subspaces by finite two-variable constraints and canonical quotient
preduals, and then identify fixed-point spaces with annihilators of
coinvariant subspaces.

\subsection{A preliminary observation}

The ambient two-Lipschitz space is already very large on infinite
metric spaces.  We record the following observation, which will also
serve to contrast the finite-dimensional example obtained later for
the positively homogeneous subspace.

\begin{proposition}\label{prop:linfty-embedding}
Let \(M\) and \(N\) be infinite pointed metric spaces.  Then the predual $\mathcal Z(M,N)=\F(M)\whot\F(N)$ of
\(\BLip_0(M,N)\) contains an isomorphic copy of \(\ell_1(\mathbb N \times \mathbb N)\).
Equivalently, \(\BLip_0(M,N)\) contains an isomorphic copy of
\(\ell_\infty(\mathbb N\times\mathbb N)\).
\end{proposition}

\begin{proof}
By \cite{OTSOLFS}, both \(\F(M)\) and \(\F(N)\) contain complemented
subspaces isomorphic to \(\ell_1\).  The projective tensor product of
complemented subspaces is complemented in the corresponding projective
tensor product; see, for instance, \cite{RAR}.  Hence
$\F(M)\whot\F(N)$
contains a complemented subspace isomorphic to $\ell_1\whot\ell_1
\cong
\ell_1(\ell_1)
\cong
\ell_1(\mathbb N\times\mathbb N)
\cong
\ell_1$.
If \(Y_0\subseteq \F(M)\whot\F(N)\) is such a complemented copy of
\(\ell_1\), the adjoint of a bounded projection onto \(Y_0\) embeds
\(Y_0^*\cong\ell_\infty\) isomorphically into
$\bigl(\F(M)\whot\F(N)\bigr)^*
\cong
\BLip_0(M,N)$.
\end{proof}

\subsection{Product actions by similarities}
\label{subsec:product-actions-similarities}
Group actions on Lipschitz-free spaces and the projections associated
with such actions have been studied in several settings; see, for
instance, \cite{PILFSIBGA,LFSATAGA}. Here we consider the tensorial
two-variable situation arising from independent similarity actions on
the two metric factors.
Let $G$ and $H$ be groups acting on the pointed metric spaces
$(M,d_M,0)$ and $(N,d_N,0)$, respectively. We assume that both
actions fix the distinguished points and act by similarities. Thus,
there are group homomorphisms
$\chi_G:G\longrightarrow(0,\infty)$ and $\chi_H:H\longrightarrow(0,\infty)$
such that
$d_M(gx,gx')=\chi_G(g)d_M(x,x')$ and
$d_N(hy,hy')=\chi_H(h)d_N(y,y')$
for all $g\in G$, $h\in H$, $x,x'\in M$, and $y,y'\in N$.

We first describe the corresponding actions on the Lipschitz-free
spaces. Fix $g\in G$ and we consider the map
$\varphi_g^M:M\longrightarrow\F(M),
~\mbox{given by}~\varphi_g^M(x)
=
\frac{1}{\chi_G(g)}\delta_{gx}$.
Since the action fixes the distinguished point,
$\varphi_g^M(0)
=
\frac{1}{\chi_G(g)}\delta_{g0}
=
\delta_0
=
0$.
Moreover, using the fact that the canonical map
$\delta_M:M\to\F(M)$ is an isometric embedding, we obtain
\[
\|\varphi_g^M(x)-\varphi_g^M(x')\|
=
\frac{1}{\chi_G(g)}
\|\delta_{gx}-\delta_{gx'}\|  =
\frac{1}{\chi_G(g)}
d_M(gx,gx')  =
d_M(x,x').
\]
Hence $\varphi_g^M$ is a base-point-preserving isometric embedding.
By the universal property of $\F(M)$, there exists a unique bounded
linear operator
$V_g^M:\F(M)\longrightarrow\F(M)$
such that
$V_g^M\delta_x
=
\varphi_g^M(x)
=
\frac{1}{\chi_G(g)}\delta_{gx},
~\mbox{for any}~ x\in M$,
and $\|V_g^M\|=\Lip(\varphi_g^M)=1$.
Similarly, for every $h\in H$, the map
$\varphi_h^N(y)
=
\frac{1}{\chi_H(h)}\delta_{hy},
~\mbox{for any}~ y\in N$,
induces a unique linear isometry
$V_h^N:\F(N)\longrightarrow\F(N)$
satisfying $V_h^N\delta_y
=
\frac{1}{\chi_H(h)}\delta_{hy}$.

The following proposition records the resulting group actions.

\begin{proposition}\label{prop:free-product-action}
The families
$(V_g^M)_{g\in G}$ and
$(V_h^N)_{h\in H}$
define representations of $G$ and $H$ by surjective linear isometries
on $\F(M)$ and $\F(N)$, respectively. Consequently, for any $(g,h)\in G\times H~$ $W_{g,h}
=
V_g^M\widehat\otimes_\pi V_h^N$,
defines a representation of $G\times H$ by surjective linear
isometries on $\mathcal Z(M,N)
=
\F(M)\widehat\otimes_\pi\F(N)$.
On elementary tensors, $W_{g,h}(\delta_x\otimes\delta_y)
=
\frac{1}{\chi_G(g)\chi_H(h)}
\delta_{gx}\otimes\delta_{hy}$.
\end{proposition}

\begin{proof}
It remains to verify the representation law and surjectivity. Let
$g_1,g_2\in G$. For every $x\in M$,
\[
V_{g_1}^M V_{g_2}^M\delta_x
=
V_{g_1}^M
\left(
\frac{1}{\chi_G(g_2)}\delta_{g_2x}
\right)=
\frac{1}{\chi_G(g_2)\chi_G(g_1)}
\delta_{g_1g_2x}=
\frac{1}{\chi_G(g_1g_2)}
\delta_{g_1g_2x}=
V_{g_1g_2}^M\delta_x,
\]
where we have used that $\chi_G$ is a group homomorphism. Since
$\operatorname{span}\{\delta_x:x\in M\}$ is dense in $\F(M)$, it
follows that
$V_{g_1}^M V_{g_2}^M=V_{g_1g_2}^M$. Also,
$V_e^M=I_{\F(M)}$.
In particular,$V_g^M V_{g^{-1}}^M
=
V_{g^{-1}}^M V_g^M
=
I_{\F(M)}$ and hence
$(V_g^M)^{-1}=V_{g^{-1}}^M$.
Thus each $V_g^M$ is a surjective linear isometry. The same argument
applies to $(V_h^N)_{h\in H}$.

Now, since the projective tensor product is functorial, the operator
$V_g^M\otimes V_h^N:
\F(M)\otimes\F(N)
\longrightarrow
\F(M)\otimes\F(N)$
extends uniquely to a bounded linear operator
$W_{g,h}
=
V_g^M\widehat\otimes_\pi V_h^N
:
\F(M)\widehat\otimes_\pi\F(N)
\longrightarrow
\F(M)\widehat\otimes_\pi\F(N)$.
Moreover,
$\|W_{g,h}\|
\leq
\|V_g^M\|\,\|V_h^N\|
=1$.
Its inverse is
$W_{g^{-1},h^{-1}}
=
V_{g^{-1}}^M\widehat\otimes_\pi V_{h^{-1}}^N,$
which is also contractive. Therefore $W_{g,h}$ is a surjective linear
isometry.

Finally, for $(g_1,h_1),(g_2,h_2)\in G\times H$,
\[
W_{g_1,h_1}W_{g_2,h_2}=
(V_{g_1}^M V_{g_2}^M)
\widehat\otimes_\pi
(V_{h_1}^N V_{h_2}^N)=
V_{g_1g_2}^M
\widehat\otimes_\pi
V_{h_1h_2}^N=
W_{g_1g_2,h_1h_2}.
\]
Thus $W:=(W_{g,h})_{(g,h)\in G\times H}$ is an action of $G\times H$
by surjective linear isometries on $\mathcal Z(M,N)$.
\end{proof}

The corresponding action on $\BLip_0(M,N)$ is defined by
\begin{equation}\label{eq:dual-product-action}
(U_{g,h}T)(x,y)
=
\chi_G(g)\chi_H(h)
T(g^{-1}x,h^{-1}y).
\end{equation}

\begin{proposition}\label{prop:BLip-product-action}
The family $U:=\{U_{g,h}:(g,h)\in G\times H\}$ is an action by surjective
linear isometries on $\BLip_0(M,N)$.  Moreover, $U_{g,h}=W_{g^{-1},h^{-1}}^*$.
\end{proposition}

\begin{proof}
Let $T\in\BLip_0(M,N)$.  Now, for $x,x'\in M$ and $y,y'\in N$,
\begin{align*}
    |U_{g,h}T(x,y)-&U_{g,h}T(x',y)
-U_{g,h}T(x,y')+U_{g,h}T(x',y')|\\
&=\chi_G(g)\chi_H(h)
\big|T(g^{-1}x,h^{-1}y)-T(g^{-1}x',h^{-1}y)-T(g^{-1}x,h^{-1}y')
+T(g^{-1}x',h^{-1}y')\big|\\
&\leq \chi_G(g)\chi_H(h)\BLip(T)
 d_M(g^{-1}x,g^{-1}x')d_N(h^{-1}y,h^{-1}y')\\
 &=\BLip(T)d_M(x,x')d_N(y,y')
\end{align*}
Thus $\BLip(U_{g,h}T)\leq\BLip(T)$.  Applying the same argument to
$(g^{-1},h^{-1})$ gives equality and surjectivity.

Finally,
$$ \left(W_{g^{-1},h^{-1}}^*T\right)(\delta_x\otimes\delta_y)=T_L\left(W_{g^{-1},h^{-1}}(\delta_x\otimes\delta_y)\right)=\chi_G(g)\chi_H(h)T(g^{-1}x,h^{-1}y),$$
which proves the adjoint relation. Recall that $T_L$ is the unique linearization map corresponding $T \in \BLip_0(M,N)$ as in \ref{eq:BLip-predual}.
\end{proof}

\subsection{Pointwise-closed invariant subspaces}

We equip $\BLip_0(M,N)$ with the topology of pointwise convergence.  A
net $(T_i)$ converges pointwise to $T$ precisely when
$T_i(x,y)\longrightarrow T(x,y)$ for any $x\in M,\ y\in N$.
Equivalently, this is the locally convex topology
$\sigma\bigl(\BLip_0(M,N),\mathscr M_2(M,N)\bigr)$.

\begin{theorem}\label{thm:pointwise-invariant-BLip}
Let $\mathcal S\subseteq\BLip_0(M,N)$ be a linear subspace which is
closed for pointwise convergence and invariant under the product action
$U$, that is,
$U_{g,h}(\mathcal S)=\mathcal S$ for any $(g,h)\in G\times H$.
Put $K=\mathcal S_\perp
=\{u\in\mathcal Z(M,N): T_L(u)=0
\text{ for every }T\in\mathcal S\}$.
Then:
\begin{enumerate}
\item $\mathcal S$ is $weak^*$ closed in
$\BLip_0(M,N)=\mathcal Z(M,N)^*$;
\item if
$\mathscr R=K\cap\mathscr M_2(M,N)$,
then
$\mathcal S=\mathscr R^\perp,~K=\overline{\mathscr R}^{\,\|\cdot\|_\pi}$;
\item $K$ is invariant under the predual action $W$;
\item there is a canonical isometric identification
$\mathcal S
\cong
\left(\bigslant{\mathcal Z(M,N)}{K}\right)^*$.
\end{enumerate}
Thus every invariant pointwise-closed linear subspace of
$\BLip_0(M,N)$ has a canonical quotient predual.
\end{theorem}

\begin{proof}
Since the topology of pointwise convergence on $\BLip_0(M,N)$ is
precisely
$\sigma\bigl(\BLip_0(M,N),\mathcal M_2(M,N)\bigr)$,
and $S$ is closed for this topology, the bipolar theorem for the dual
pair
$\bigl(\BLip_0(M,N),\mathcal M_2(M,N)\bigr)$
gives
$S=(S^{\perp,\mathcal M_2})^\perp$,
where
$S^{\perp,\mathcal M_2}
=
\{u\in\mathcal M_2(M,N):T_L(u)=0\ \text{for every }T\in S\}$.

Since
$\mathcal S_{\perp,\mathscr M_2}
=
\mathcal S_\perp\cap\mathscr M_2(M,N)$,
we obtain
$\mathcal S
=
\bigl(
\mathcal S_\perp\cap\mathscr M_2(M,N)
\bigr)^\perp
=
\mathscr R^\perp$.
Now, taking preannihilators and using the bipolar theorem we have
$K=\mathcal S_\perp
=\overline{\mathscr R}^{\,\|\cdot\|_\pi}$.

Let $u\in K$, $T\in\mathcal S$, and $(g,h)\in G\times H$.  Then
$ T_L(W_{g,h}u)
=W_{g,h}^*T(u)
= U_{g^{-1},h^{-1}}T(u)
=0$,
because $\mathcal S$ is invariant under $U$.  Hence $W_{g,h}u\in K$.
Finally, $\mathcal S=K^\perp$. Therefore, from the standard duality for quotients we get
$K^\perp\cong(\mathcal Z(M,N)/K)^*$.
\end{proof}

\begin{remark}\label{rem:finite-constraints}
The theorem gives a concrete description by finite two-variable
constraints.  Every $u\in\mathscr R$ has the form
$u=\sum_{j=1}^n a_j\delta_{x_j}\otimes\delta_{y_j}$,
and the equation $T_L(u)=0$ is exactly
$\sum_{j=1}^n a_jT(x_j,y_j)=0$.
Thus the pointwise-closed subspace is the solution space of a family of
finite constraints.  Invariance of $\mathcal S$ is reflected by
invariance of the norm-closed relation space $K$ under the predual
action $W$.  This is the two-variable tensor counterpart of the
finite-constraint/quotient-predual viewpoint in \cite{DouchaIPC}.
\end{remark}


\subsection{Fixed points and the coinvariant tensor quotient}

The quotient construction below is the Banach-space counterpart of the
standard duality between invariants of a contragredient representation
and coinvariants. Recall that, for a representation of a group $G$ on a
vector space $X$, the coinvariant space is obtained by quotienting $X$
by the linear span of the differences $gx-x$. In the algebraic setting,
the dual of the coinvariant space identifies naturally with the invariant
subspace of the contragredient representation; see \cite{SerreRep},
\cite[Section~2.2.2 and Exercise~2.2.19]{EK}.

In the present setting, the representation acts on the Banach space
$\mathcal Z(M,N)=\F(M)\widehat{\otimes}_{\pi}\F(N)$,
and we therefore pass to the norm-closed coinvariant subspace
\[
\Coinv(W)
=
\overline{\operatorname{span}}
\left\{
W_{g,h}u-u:
(g,h)\in G\times H,\ u\in\mathcal Z(M,N)
\right\}.
\]
Theorem~\ref{thm:fixed-point-coinvariant} is the corresponding
Banach-space formulation for the tensor representation arising from
two-Lipschitz mappings.

\begin{theorem}[Fixed-point/coinvariant duality for the tensor representation]
\label{thm:fixed-point-coinvariant}
The fixed-point space $\Fix(U)$ is $weak^*$ closed and
$\Fix(U)=\Coinv(W)^\perp$.
Consequently,
\begin{equation}\label{eq:fixed-coinvariant-duality}
\Fix(U)
\cong
\left(
\bigslant{\mathcal Z(M,N)}{\Coinv(W)}
\right)^*
\end{equation}
isometrically.  Hence the coinvariant tensor quotient is a canonical
predual of the fixed-point space. Equivalently, it is a canonical predual of $\BLip_0^{G\times H,\chi}(M,N)$.
\end{theorem}

\begin{proof}
For $T\in\BLip_0(M,N)$, the condition $T\in\Fix(U)$ is equivalent,
after replacing $(g,h)$ by $(g^{-1},h^{-1})$, to
$W_{g,h}^*T_L=T_L
~\mbox{for any}~ (g,h)\in G\times H$.
This is equivalent to
$ T_L\left(W_{g,h}u-u\right)=0$
for every $(g,h)\in G\times H$ and every $u\in\mathcal Z(M,N)$.
Therefore,
$\Fix(U)=\Coinv(W)^\perp$.
Now, an annihilator is $weak^*$ closed. Thus from the quotient-duality formula
we have \eqref{eq:fixed-coinvariant-duality}.
\end{proof}
Thus the standard invariant-coinvariant duality, when applied to the
tensor representation $W$, yields a canonical quotient predual for the
fixed-point subspace of $\BLip_0(M,N)$. The specific content in the
present setting is that this fixed-point space is precisely the \emph{space of
character-equivariant two-Lipschitz mappings}; see
Remark~\ref{rem:fixed-character-equivariant}.

\begin{corollary}\label{cor:fixed-point-pointwise-closed}
The fixed-point space $\Fix(U)$ is closed for pointwise convergence.
\end{corollary}

\begin{proof}
For each $(g,h)\in G\times H$, the equation $U_{g,h}T=T$ is a family of
pointwise equations.  These equations are preserved under pointwise
limits, and their intersection is $\Fix(U)$. Hence the proof follows.
\end{proof}
\begin{remark}[Character-equivariant interpretation]
\label{rem:fixed-character-equivariant}
Recall that the product character associated with the similarity
actions of $G$ and $H$ is
$\chi(g,h)=\chi_G(g)\chi_H(h)$, for any $ (g,h)\in G\times H,$
and that the induced action on $\BLip_0(M,N)$ is given by
$(U_{g,h}T)(x,y)
=
\chi(g,h)\,
T(g^{-1}x,h^{-1}y)$.
Thus, for $T\in\BLip_0(M,N)$,
 we have $U_{g,h}T=T$ for every $(g,h)\in G\times H$
if and only if
$T(gx,hy)
=
\chi(g,h)T(x,y)$
for every $g\in G$, $h\in H$, $x\in M$, and $y\in N$.
Consequently,
$\Fix(U)
=
\BLip_0^{G\times H,\chi}(M,N)$.

Hence Theorem~\ref{thm:fixed-point-coinvariant} may equivalently be written as
$\BLip_0^{G\times H,\chi}(M,N)
=
\Coinv(W)^\perp$,
and therefore
$\BLip_0^{G\times H,\chi}(M,N)
\cong
\left(
\bigslant{\mathcal Z(M,N)}{\Coinv(W)}
\right)^*$. In particular, the character-equivariant two-Lipschitz
space is $weak^*$ closed in $\BLip_0(M,N)$ and the coinvariant tensor
quotient
$\bigslant{\mathcal Z(M,N)}{\Coinv(W)}$
is a canonical predual of it.

Moreover, Corollary~\ref{cor:fixed-point-pointwise-closed} shows that
$\BLip_0^{G\times H,\chi}(M,N)$ is also closed for pointwise
convergence. Thus the character-equivariant class provides a natural
example of an invariant pointwise-closed, and hence $weak^*$ closed,
subspace of the two-Lipschitz space.
\end{remark}

\subsection{Amenability and dual-valued contractive projections}

We now pass to a Banach-valued version of the averaging argument.
This is most convenient when the codomain is a dual space. The groups
below are regarded as discrete groups. Recall that a discrete group
$\Gamma$ is amenable if there exists a left-invariant mean
$m:\ell_\infty(\Gamma)\longrightarrow \mathbb R$,
that is, a positive linear functional of norm one satisfying
$m(\mathbf 1)=1
$ and
$m(L_\gamma\varphi)=m(\varphi)
~~\mbox{for any}~\gamma\in\Gamma,\ \varphi\in\ell_\infty(\Gamma)$,
where $(L_\gamma\varphi)(\eta)=\varphi(\gamma^{-1}\eta)$.
Since the product of two amenable groups is amenable, if $G$ and $H$
are amenable then there exists a left-invariant mean on
$\ell_\infty(G\times H)$.

For an arbitrary Banach space \(F\) we define
\[
\Fix(U;F)
=
\{T\in\BLip_0(M,N;F):U_{g,h}T=T
\text{ for all }(g,h)\in G\times H\}.
\]
\begin{theorem}\label{thm:amenable-projection}
Let $G$ and $H$ be amenable groups, and let $E$ be a Banach space.
Then there exists a contractive linear projection
$P_m:
\BLip_0(M,N;E^*)
\longrightarrow
\Fix(U;E^*)$.
More precisely, if
$m:\ell_\infty(G\times H)\longrightarrow\mathbb R$
is a left-invariant mean, then $P_mT(x,y)\in E^*$ is determined by
\begin{equation}\label{eq:invariant-mean-projection}
\bigl(P_mT(x,y)\bigr)(e)
=
m\Bigl(
(g,h)\longmapsto
\bigl((U_{g,h}T)(x,y)\bigr)(e)
\Bigr),~~\mbox{for any}~ e\in E.
\end{equation}
\end{theorem}

\begin{proof}
Since $G$ and $H$ are amenable, $G\times H$ is amenable.  Let $m$ be a
left-invariant mean.  For $T\in\BLip_0(M,N;E^*)$, from the vanishing of $T$
on the coordinate axes we get
\[
\|T(u,v)\|
\leq
\BLip(T)d_M(u,0)d_N(v,0).
\]
Hence, for $e\in E$, we have
$$| \left((U_{g,h}T)(x,y)\right)(e)|\leq \|(U_{g,h}T)(x,y)\|\,\|e\|\leq \BLip(T)d_M(x,0)d_N(y,0)\|e\|.$$
Thus for fixed $e \in E$ the function $(g,h)\longmapsto
\bigl((U_{g,h}T)(x,y)\bigr)(e)$ is bounded.  Now, for
fixed $(x,y)$, its mean is a bounded linear functional of $e$, so it
defines an element $P_mT(x,y)\in E^*$.

Let $x,x'\in M$, $y,y'\in N$, and $e\in B_E$.  By positivity and norm
one of the invariant mean, we get
\begin{align*}
&|\left( P_mT(x,y)-P_mT(x',y)-P_mT(x,y')+P_mT(x',y')\right)(e)| \\
&\leq m\Bigl((g,h)\longmapsto
\bigl\|(U_{g,h}T)(x,y)-(U_{g,h}T)(x',y)-(U_{g,h}T)(x,y')+(U_{g,h}T)(x',y')\bigr\|\Bigr)\\
&\leq \BLip(T)d_M(x,x')d_N(y,y').
\end{align*}
Taking the supremum over $e\in B_E$ we have
$\BLip(P_mT)\leq\BLip(T),$
so $\|P_m\|\leq1$.

Let $(k,l)\in G\times H$.  Since $U$ is a representation,
$U_{k,l}U_{g,h}=U_{kg,lh}$.
Now, applying left invariance of $m$ we have
$U_{k,l}P_mT=P_mT$.
Thus $P_mT\in\Fix(U;E^*)$.  

If $T\in\Fix(U;E^*)$, then the function being averaged in \eqref{eq:invariant-mean-projection} is constant, and hence
$P_mT=T$.  Therefore, $P_m^2=P_m$ and
$\operatorname{ran}P_m=\Fix(U;E^*)$.  If this range is nonzero, every
projection onto it has norm at least one, so $\|P_m\|=1$.
\end{proof}

\begin{corollary}\label{cor:amenable-scalar-projection}
If $G$ and $H$ are amenable, then there exists a contractive projection
\[
P_m:\BLip_0(M,N)\longrightarrow\Fix(U).
\]
If $\Fix(U)\neq\{0\}$, then $\|P_m\|=1$.
\end{corollary}

\begin{remark}\label{rem:mean-not-wstar}
The projection supplied by Theorem~\ref{thm:amenable-projection}
is an invariant-mean projection.  It is generally nonconstructive and
depends on the choice of the mean.  In the scalar case, it should not in
general be expected to be $weak^*$ continuous; in particular, the
construction does not automatically provide a bounded projection on the
coinvariant predual whose adjoint is $P_m$.
\end{remark}

\begin{remark}\label{rem:general-codomain-bidual}
Let $F$ be an arbitrary Banach space, and assume that $G$ and $H$ are
amenable. We fix a left-invariant mean $m$ on $\ell_\infty(G\times H)$.
For $T\in\BLip_0(M,N;F)$, we again recall that the induced action is given by
$(U_{g,h}T)(x,y)
=
\chi_G(g)\chi_H(h)
T(g^{-1}x,h^{-1}y), ~\mbox{for any}~ (g,h)\in G\times H$. For each $(x,y)\in M\times N$, we define a functional
$\Phi_{T,x,y}:F^*\longrightarrow\mathbb R$
by
$\Phi_{T,x,y}(f^*)
=
m\Bigl(
(g,h)\longmapsto
f^{*}\bigl((U_{g,h}T)(x,y)\bigr)
\Bigr), ~~\mbox{for any}~ f^*\in F^*$.

We first verify that $\Phi_{T,x,y}\in F^{**}$. Since
$T$ vanishes on the coordinate axes,
$\|T(u,v)\|
\leq
\BLip(T)\,d_M(u,0)d_N(v,0),
~~\mbox{for any}~ u\in M,\ v\in N$.
As the actions fix the distinguished points and have similarity
characters $\chi_G$ and $\chi_H$, we have
\[
d_M(g^{-1}x,0)
=
\chi_G(g)^{-1}d_M(x,0),
\qquad
d_N(h^{-1}y,0)
=
\chi_H(h)^{-1}d_N(y,0).
\]
Consequently,
$\|(U_{g,h}T)(x,y)\|
=
\chi_G(g)\chi_H(h)
\|T(g^{-1}x,h^{-1}y)\|\leq
\BLip(T)d_M(x,0)d_N(y,0)$.
Hence, for every $f^*\in F^*$,
$\left|
f^{*}\bigl((U_{g,h}T)(x,y)\bigr)
\Bigr)\right|
\leq
\BLip(T)d_M(x,0)d_N(y,0)\|f^*\|$.
Thus the scalar-valued function
$(g,h)\longmapsto
f^{*}\bigl((U_{g,h}T)(x,y)\bigr)
\Bigr)$
belongs to $\ell_\infty(G\times H)$. Since $m$ has norm one, we have
$|\Phi_{T,x,y}(f^*)|
\leq
\BLip(T)d_M(x,0)d_N(y,0)\|f^*\|$.
Moreover, $\Phi_{T,x,y}$ is linear in $f^*$. Therefore
$\Phi_{T,x,y}$ is a bounded linear functional on $F^*$, with
$\|\Phi_{T,x,y}\|
\leq
\BLip(T)d_M(x,0)d_N(y,0)$.
Hence
$\Phi_{T,x,y}\in(F^*)^*=F^{**}$.

We may therefore define
$(\widetilde P_m^F T)(x,y):=\Phi_{T,x,y},
~\mbox{for any} (x,y)\in M\times N$.
Thus
$\widetilde P_m^F T:M\times N\longrightarrow F^{**}$.
Further, we can deduce that 
$\BLip(\widetilde P_m^F T)\leq\BLip(T)$,
and hence
$\widetilde P_m^F:
\BLip_0(M,N;F)
\longrightarrow
\BLip_0(M,N;F^{**})$
is a contractive linear map. By the left invariance of $m$, its image
consists of fixed points for the induced action on
$\BLip_0(M,N;F^{**})$.

If $J_F:F\to F^{**}$ denotes the canonical embedding, then for every
$T \in \Fix(U;F)$ we get,
$\widetilde P_m^F T=J_F\circ T$
In particular, if $F$ is reflexive, then $J_F$ is surjective and
$P_m^F:=J_F^{-1}\circ\widetilde P_m^F$
defines a contractive projection of $\BLip_0(M,N;F)$ onto its
fixed-point subspace.
\end{remark}

\subsection{The dilation action and positive homogeneity}

 We now consider the groups $G=H=(0,\infty)$
with multiplication and let these groups act by dilations on the Banach spaces $X,~Y$, respectively,
$a\cdot x=ax$, and $b\cdot y=by$.
Then $\chi_G(a)=a$ and $\chi_H(b)=b$.  Also, the induced action on
$\BLip_0(X,Y)$ is
\begin{equation}\label{eq:dilation-action}
(U_{a,b}T)(x,y)
=
ab\,T\left(\frac{x}{a},\frac{y}{b}\right).
\end{equation}

Further, we recall that $\BLip_0^{ph}(X,Y)$ consists of those
$T\in\BLip_0(X,Y)$ satisfying
$T(ax,y)=aT(x,y)$ and $T(x,by)=bT(x,y)$
for every $a,b\geq0$.

\begin{proposition}\label{prop:ph-fixed}
We have
\[
\BLip_0^{ph}(X,Y)
=
\Fix\bigl((0,\infty)^2,\BLip_0(X,Y)\bigr).
\]
In particular, $\BLip_0^{ph}(X,Y)$ is pointwise closed in
$\BLip_0(X,Y)$.
\end{proposition}

\begin{proof}
Let $T$ is fixed by the action, then
$ab\,T\left(\frac{x}{a},\frac{y}{b}\right)=T(x,y)$ for any $(a,b>0)$.
Replacing $x$ by $ax$ and $y$ by $by$ we have
$T(ax,by)=abT(x,y)$.
Taking $b=1$ and $a=1$, respectively, yields separate positive
homogeneity.  The case of scalar $0$ follows from the vanishing of $T$
on the coordinate axes.  Conversely, separate positive homogeneity
immediately implies \eqref{eq:dilation-action} fixes $T$.

The pointwise-closedness follows either directly by passing to
pointwise limits in the homogeneity identities or from
Corollary~\ref{cor:fixed-point-pointwise-closed}.
\end{proof}
\section{Applications}
\medskip
\noindent\textbf{The predual dilation action.}
According to the general construction of the predual action as in the Proposition \ref{prop:free-product-action} and its preceding discussion, we have
$V_a^X\delta_x=\frac1a\delta_{ax},
~V_b^Y\delta_y=\frac1b\delta_{by}$
and
$W_{a,b}
=
V_a^X\widehat\otimes_\pi V_b^Y$,
defines a representation of $(0,\infty)^2$ by surjective linear
isometries on
$\mathcal Z(X,Y)
=
\F(X)\widehat\otimes_\pi\F(Y)$.
In particular,
$W_{a,b}(\delta_x\otimes\delta_y)
=
\frac1{ab}\,
\delta_{ax}\otimes\delta_{by}$.

Again, we recall that the associated coinvariant subspace is
$$\Coinv(W)
=
\overline{\operatorname{span}}
\left\{
W_{a,b}u-u:
a,b>0,\;
u\in\F(X)\widehat\otimes_\pi\F(Y)
\right\}.$$

We put $N_{ph}(X)
=
\overline{\operatorname{span}}
\left\{
\delta_{ax}-a\delta_x:
a>0,\ x\in X
\right\}
\subseteq\F(X)$ and
$N_{ph}(Y)
=
\overline{\operatorname{span}}
\left\{
\delta_{by}-b\delta_y:
b>0,\ y\in Y
\right\}
\subseteq\F(Y)$.
Further, we consider the closed relation space
\begin{equation}\label{eq:dilation-relation-space}
K_{ph}(X,Y)
=
\overline{
N_{ph}(X)\otimes\F(Y)
+
\F(X)\otimes N_{ph}(Y)
}^{\,\|\cdot\|_\pi}
\subseteq
\F(X)\widehat\otimes_\pi\F(Y).
\end{equation}
Here the tensor products inside the closure are algebraic tensor
products viewed as subspaces of the ambient projective tensor product.

\begin{proposition}\label{prop:dilation-coinvariant}
For the product dilation action of $(0,\infty)^2$, the general
coinvariant subspace coincides with the relation space $K_{ph}(X,Y)$: $~\Coinv(W)=K_{ph}(X,Y)$.
Equivalently,
\begin{equation}\label{eq:H-character-relations}
K_{ph}(X,Y)
=
\overline{\operatorname{span}}
\left\{
\delta_{ax}\otimes\delta_{by}
-
ab\,\delta_x\otimes\delta_y:
a,b>0,\ x\in X,\ y\in Y
\right\}^{\,\|\cdot\|_\pi}.
\end{equation}

Moreover,
$\BLip_0^{ph}(X,Y)
=
\BLip_0^{(0,\infty)^2,\chi}(X,Y)
=
\Fix(U)
=
K_{ph}(X,Y)^\perp,$
where $\chi(a,b)=ab$.
Consequently,
\begin{equation}\label{eq:ph-quotient-predual}
\BLip_0^{ph}(X,Y)
\cong
\left(
\bigslant{\F(X)\widehat\otimes_\pi\F(Y)}{K_{ph}(X,Y)}
\right)^*
\end{equation}
isometrically. In particular, $\bigslant{\F(X)\widehat\otimes_\pi\F(Y)}{K_{ph}(X,Y)}$
is the canonical coinvariant predual of the positively homogeneous
two-Lipschitz space.
\end{proposition}

\begin{proof}
We first identify the coinvariant relation space. For $a>0$ and
$x\in X$,
\[
(V_a^X-I)\delta_x
=
\frac1a\delta_{ax}-\delta_x
=
\frac1a(\delta_{ax}-a\delta_x).
\]
Since $\operatorname{span}\{\delta_x:x\in X\}$
is dense in $\F(X)$, it follows that
$N_{ph}(X)
=
\overline{\operatorname{span}}
\left\{
(V_a^X-I)\mu:
a>0,\ \mu\in\F(X)
\right\}$.
Similarly, $N_{ph}(Y)
=
\overline{\operatorname{span}}
\left\{
(V_b^Y-I)\nu:
b>0,\ \nu\in\F(Y)
\right\}$.

Now, for $a,b>0$, we have
\begin{equation}\label{eq:dilation-W-decomposition}
W_{a,b}-I
=
(V_a^X-I)\widehat\otimes_\pi V_b^Y
+
I\widehat\otimes_\pi(V_b^Y-I).
\end{equation}
Indeed, on an elementary tensor $\mu\otimes\nu$,
$$
(W_{a,b}-I)(\mu\otimes\nu)=
V_a^X\mu\otimes V_b^Y\nu-\mu\otimes\nu =
(V_{a}^{X}-I)\mu\otimes V_b^Y\nu
+
\mu\otimes(V_b^Y-I)\nu.$$
The equality therefore extends, by continuity, to
$\F(X)\widehat\otimes_\pi\F(Y)$.

The first term on the right-hand side of
\eqref{eq:dilation-W-decomposition} belongs to the norm closure of
$N_{ph}(X)\otimes\F(Y)$,
whereas the second belongs to the norm closure of
$\F(X)\otimes N_{ph}(Y)$.
Hence
$(W_{a,b}-I)
\bigl(\F(X)\widehat\otimes_\pi\F(Y)\bigr)
\subseteq K_{ph}(X,Y)$,
and therefore $\Coinv(W)\subseteq H$.

For the reverse inclusion, first we observe that
$W_{a,1}-I
=
(V_a^X-I)\widehat\otimes_\pi I$
and $W_{1,b}-I
=
I\widehat\otimes_\pi(V_b^Y-I)$.
Thus, for $\mu\in\F(X)$ and $\nu\in\F(Y)$,
\[
(V_a^X-I)\mu\otimes\nu
=
(W_{a,1}-I)(\mu\otimes\nu)
\in\Coinv(W),
~\mbox{and}~\mu\otimes(V_b^Y-I)\nu
=
(W_{1,b}-I)(\mu\otimes\nu)
\in\Coinv(W).
\]
Using the descriptions of $N_{ph}(X)$ and $N_{ph}(Y)$ above, together with the
fact that $\Coinv(W)$ is norm closed, we obtain
$N_{ph}(X)\otimes\F(Y)
+
\F(X)\otimes N_{ph}(Y)
\subseteq
\Coinv(W)$
after taking the appropriate norm closure. Hence $K_{ph}(X,Y)\subseteq\Coinv(W)$. Therefore
$\Coinv(W)=K_{ph}(X,Y)$.

We next relate this equality to the general character-equivariant
theory. On canonical elementary tensors,
\[
(W_{a,b}-I)(\delta_x\otimes\delta_y)
=
\frac1{ab}\delta_{ax}\otimes\delta_{by}
-
\delta_x\otimes\delta_y=
\frac1{ab}
\left(
\delta_{ax}\otimes\delta_{by}
-
ab\,\delta_x\otimes\delta_y
\right).
\]
Since the span of the elementary tensors
$\delta_x\otimes\delta_y$ is dense in
$\F(X)\widehat\otimes_\pi\F(Y)$, this proves
\eqref{eq:H-character-relations}.

By the general fixed-point/coinvariant duality theorem, $\Fix(U)
=
\Coinv(W)^\perp$.
Since $\Coinv(W)=H$, we obtain
$\Fix(U)=K_{ph}(X,Y)^\perp$.

On the other hand, for the present dilation action the product
character is $\chi(a,b)=ab$,
and the induced action on $\BLip_0(X,Y)$ is $(U_{a,b}T)(x,y)
=
ab\,T\left(\frac{x}{a},\frac{y}{b}\right)$.
Hence, by the general character-equivariant interpretation of the
fixed-point space,
$\Fix(U)
=
\BLip_0^{(0,\infty)^2,\chi}(X,Y)$.
Finally, the $\chi$-equivariance relation
$T(ax,by)=ab\,T(x,y)$
is equivalent to positive homogeneity in each variable. Therefore,
$\BLip_0^{(0,\infty)^2,\chi}(X,Y)
=
\BLip_0^{ph}(X,Y)$.
Combining the preceding identities gives
$\BLip_0^{ph}(X,Y)
=
\BLip_0^{(0,\infty)^2,\chi}(X,Y)
=
\Fix(U)
=
K_{ph}(X,Y)^\perp$.

The quotient-duality part of the general
fixed-point/coinvariant theorem now yields
$K_{ph}(X,Y)^\perp
\cong
\left(
\bigslant{\F(X)\widehat\otimes_\pi\F(Y)}{K_{ph}(X,Y)}
\right)^*$, which proves
\eqref{eq:ph-quotient-predual}.
\end{proof}

\begin{theorem}\label{thm:ph-vector-valued-quotient}
Let \(E\) be a Banach space.  Linearization induces a canonical
isometric isomorphism
\begin{equation}\label{eq:ph-vector-valued-quotient}
\BLip_0^{ph}(X,Y;E)
\cong
\Lop\left(
\bigslant{\F(X)\whot\F(Y)}{K_{ph}(X,Y)},E
\right).
\end{equation}
Equivalently, a map \(T\in\BLip_0(X,Y;E)\) is positively homogeneous
if and only if its linearization
$T_L:\F(X)\whot\F(Y)\longrightarrow E$
vanishes on \(K_{ph}(X,Y)\).
\end{theorem}

\begin{proof}
Let \(T_L\) denote the linearization of
\(T\in\BLip_0(X,Y;E)\).  Now, positive homogeneity in the first variable is
equivalent to
$T_L\bigl((\delta_{ax}-a\delta_x)\otimes\delta_y\bigr)=0$
for every \(a\geq0\), \(x\in X\), and \(y\in Y\).  Similarly, positive
homogeneity in the second variable is equivalent to
$T_L\bigl(\delta_x\otimes(\delta_{by}-b\delta_y)\bigr)=0$
for every \(b\geq0\), \(x\in X\), and \(y\in Y\).
By linearity, density, and continuity, these two families of relations
are equivalent to
$T_L|_H=0$.
Hence \(T_L\) factors uniquely through the quotient map
\[
q:\F(X)\whot\F(Y)\longrightarrow
\bigl(\F(X)\whot\F(Y)\bigr)/K_{ph}(X,Y).
\]
The standard quotient factorization theorem gives the isometric
identification \eqref{eq:ph-vector-valued-quotient}.
\end{proof}

We further put
$\F^{ph}(X)=\F(X)/N_{ph}(X), ~ \mbox{and}~\F^{ph}(Y)=\F(Y)/N_{ph}(Y)$. These are the one-variable positively homogeneous quotient
preduals considered in \cite{AM2}.

\begin{corollary}\label{cor:ph-tensor-predual}
There is a canonical isometric isomorphism
\begin{equation}\label{eq:ph-tensor-quotient}
\bigslant{\F(X)\whot\F(Y)}{K_{ph}(X,Y)}
\cong
\F^{ph}(X)\whot\F^{ph}(Y).
\end{equation}
Consequently, for every Banach space \(E\),
\begin{equation}\label{eq:ph-vector-tensor-linearization}
\BLip_0^{ph}(X,Y;E)
\cong
\Lop\bigl(\F^{ph}(X)\whot\F^{ph}(Y),E\bigr).
\end{equation}
In particular,
\begin{equation}\label{eq:ph-tensor-dual}
\BLip_0^{ph}(X,Y)
\cong
\bigl(\F^{ph}(X)\whot\F^{ph}(Y)\bigr)^*,
\end{equation}
so \(\BLip_0^{ph}(X,Y)\) is a $weak^*$ closed subspace
of \(\BLip_0(X,Y)\) with canonical predual
\(\F^{ph}(X)\whot\F^{ph}(Y)\).
\end{corollary}

\begin{proof}
Apply the quotient property of the projective tensor product to the
quotient maps
\[
q_X:\F(X)\to\F(X)/N_{ph}(X),
\qquad
q_Y:\F(Y)\to\F(Y)/N_{ph}(Y).
\]
The kernel of \(q_X\whot q_Y\) is precisely the closed subspace \(K_{ph}(X,Y)\)
in \eqref{eq:dilation-relation-space}.  This gives
\eqref{eq:ph-tensor-quotient}.  Formula
\eqref{eq:ph-vector-tensor-linearization} follows from
Theorem~\ref{thm:ph-vector-valued-quotient}, and the scalar statement is
its special case \(E=\mathbb R\).
\end{proof}

\subsection{Restriction to the unit spheres and the radial decomposition}

For a Banach space \(X\), we put
\[
\mathbb{S}_X^0=\mathbb{S}_X\cup\{0\},
\qquad
\sigma_X(x)=
\begin{cases}
x/\|x\|,&x\neq0,\\
0,&x=0.
\end{cases}
\]
The following elementary estimate will be used repeatedly.

\begin{lemma}[Sphere estimate]\label{lem:sphere-estimate}
For all \(x,x'\in X\), it follows $\|x\|\,
\|\sigma_X(x)-\sigma_X(x')\|
\leq
2\|x-x'\|$.
\end{lemma}

\begin{proof}
If \(x=0\), the assertion is immediate.  If \(x\neq0\) and \(x'=0\),
the left-hand side equals \(\|x\|\).  Finally, if \(x,x'\neq0\), then
$$\|x\|\,\|\sigma_X(x)-\sigma_X(x')\|
=
\left\|x-\frac{\|x\|}{\|x'\|}x'\right\| \leq
\|x-x'\|+
\left|\,\|x'\|-\|x\|\,\right| \leq
2\|x-x'\|.$$
\end{proof}

\begin{theorem}[Sphere restriction]\label{thm:sphere-restriction}
Let \(X,Y,E\) be Banach spaces.  The restriction map
\[
\Lambda:
\BLip_0^{ph}(X,Y;E)
\longrightarrow
\BLip_0(\mathbb{S}_X^0,\mathbb{S}_Y^0;E),
~ \mbox{given by}~\Lambda(T)=T|_{\mathbb{S}_X^0\times \mathbb{S}_Y^0},
\]
is a linear homeomorphism.  More precisely,
$\BLip(\Lambda T)\leq\BLip(T)$
and its inverse is the radial extension
\begin{equation}\label{eq:radial-extension}
(\mathcal R G)(x,y)
=
\|x\|\|y\|\,G(\sigma_X(x),\sigma_Y(y)), ~~~G \in \BLip_0(\mathbb{S}_X^0,\mathbb{S}_Y^0;E).
\end{equation}
for which
$\BLip(\mathcal R G)\leq9\,\BLip(G)$.
Consequently,
$\BLip(\Lambda T)
\leq
\BLip(T)
\leq
9\,\BLip(\Lambda T)$ for all $T\in\BLip_0^{ph}(X,Y;E)$.
\end{theorem}

\begin{proof}
The restriction map is clearly linear and contractive.  If
\(T\in\BLip_0^{ph}(X,Y;E)\), then
\[
T(x,y)
=
\|x\|\|y\|\,
T(\sigma_X(x),\sigma_Y(y)),
\]
so \(\Lambda\) is injective.

Let \(G\in\BLip_0(\mathbb{S}_X^0,\mathbb{S}_Y^0;E)\) and define \(\mathcal R G\) by
\eqref{eq:radial-extension}.  It is immediate that \(\mathcal R G\)
vanishes on the coordinate axes and is positively homogeneous in each
variable.  We put
$p=\|x\|,~ p'=\|x'\|,~
q=\|y\|,~ q'=\|y'\|$,
and write $\widehat x=\sigma_X(x),~\widehat x'=\sigma_X(x'),~
\widehat y=\sigma_Y(y),~
\widehat y'=\sigma_Y(y')$.
By a direct regrouping we have
\begin{align*}
\mathcal RG(x,y)-\mathcal RG(x',y)
&-\mathcal RG(x,y')+\mathcal RG(x',y')\\
&=pq\bigl[
G(\widehat x,\widehat y)-G(\widehat x',\widehat y)
-G(\widehat x,\widehat y')+G(\widehat x',\widehat y')
\bigr]\\
&+(p-p')q\bigl[
G(\widehat x',\widehat y)-G(\widehat x',\widehat y')
\bigr]+p(q-q')\bigl[
G(\widehat x,\widehat y')-G(\widehat x',\widehat y')
\bigr]\\
&+(p-p')(q-q')G(\widehat x',\widehat y').
\end{align*}
Since $G$ vanishes on the coordinate axes, fixing one variable gives
a Lipschitz map in the other variable with Lipschitz constant at most
$\BLip(G)$. Moreover, for every
$u\in \mathbb{S}_X^0$ and $v\in \mathbb{S}_Y^0$,
\[
\|G(u,v)\|
\leq
\BLip(G)\|u\|\|v\|
\leq
\BLip(G).
\]
Now, using Lemma~\ref{lem:sphere-estimate} in both variables and the reverse
triangle inequality, the four terms above are bounded respectively by
$$4\BLip(G)\|x-x'\|\|y-y'\|,~2\BLip(G)\|x-x'\|\|y-y'\|,~
2\BLip(G)\|x-x'\|\|y-y'\|,
~\BLip(G)\|x-x'\|\|y-y'\|.$$
Further, their sum is at most
$9\BLip(G)\|x-x'\|\|y-y'\|$.
Thus \(\mathcal RG\in\BLip_0^{ph}(X,Y;E)\) and
\(\BLip(\mathcal RG)\leq9\BLip(G)\).  Finally,
\(\Lambda(\mathcal RG)=G\), which proves the result.
\end{proof}

\begin{proposition}[Radial projection]\label{prop:radial-projection}
Let \(X,Y,E\) be Banach spaces. 
We consider the restriction map $\operatorname{Res}_{\mathbb{S}_X^0\times \mathbb{S}_Y^0}:
\BLip_0(X,Y;E)
\longrightarrow
\BLip_0(\mathbb{S}_X^0,\mathbb{S}_Y^0;E)$ given by $T \mapsto T|_{\mathbb{S}_X^0\times \mathbb{S}_Y^0}$.
Then the operator
\[
P_{\mathrm{rad}}
=
\mathcal R\circ\operatorname{Res}_{\mathbb{S}_X^0\times \mathbb{S}_Y^0}
:
\BLip_0(X,Y;E)
\longrightarrow
\BLip_0^{ph}(X,Y;E)
\]
is a bounded linear projection with
$\|P_{\mathrm{rad}}\|\leq9$.
Moreover,
$\ker P_{\mathrm{rad}}
=
\{T\in\BLip_0(X,Y;E):
T|_{\mathbb{S}_X\times \mathbb{S}_Y}=0\}$,
and hence
$\BLip_0(X,Y;E)
=
\ker P_{\mathrm{rad}}
\oplus
\BLip_0^{ph}(X,Y;E)$
as a topological direct sum.
\end{proposition}

\begin{proof}
The restriction operator has norm at most one and
Theorem~\ref{thm:sphere-restriction} gives
\(\|\mathcal R\|\leq9\).  Hence
\(\|P_{\mathrm{rad}}\|\leq9\).  If
\(T\in\BLip_0^{ph}(X,Y;E)\), positive homogeneity gives
\(P_{\mathrm{rad}}T=T\), so \(P_{\mathrm{rad}}\) is a projection onto
\(\BLip_0^{ph}(X,Y;E)\).

By definition, \(P_{\mathrm{rad}}T=0\) if and only if the restriction of
\(T\) to \(\mathbb{S}_X\times \mathbb{S}_Y\) is zero.  The direct-sum decomposition is the
standard decomposition associated with a bounded projection.
\end{proof}

\subsection{Norm-one projections onto the positively homogeneous part}

We write $\BLip_0^{ph}(X,Y;E)$ for the
subspace of $\BLip_0(X,Y;E)$ consisting of maps that are separately
positively homogeneous in the two variables.

The product dilation group $(0,\infty)^2$ is abelian and hence
amenable.  The general averaging theorem therefore gives the following
strengthening of the scalar-valued projection result.

\begin{theorem}\label{thm:amenable-projection}
Let $X$ and $Y$ be nonzero Banach spaces and let $E$ be a Banach space.
There exists a contractive linear projection
\[
P_m:\BLip_0(X,Y;E^*)
\longrightarrow
\BLip_0^{ph}(X,Y;E^*).
\]
If $m$ is a left-invariant mean on $\ell_\infty((0,\infty)^2)$, then
$P_mT(x,y)\in E^*$ is determined by
\begin{equation}\label{eq:dual-valued-dilation-mean}
\left( P_mT(x,y)\right)(e)=
m_{(a,b)}\left[(a,b) \longmapsto
ab T\left(\frac{x}{a},\frac{y}{b}\right)(e)\right], ~ \mbox{for any}~ e\in E.
\end{equation}

Moreover,
$\BLip(P_mT)\leq\BLip(T)$.
If $E\neq\{0\}$, then $\|P_m\|=1$.
\end{theorem}

\begin{proof}
Apply Theorem~\ref{thm:amenable-projection} to the dilation
action.  Proposition~\ref{prop:ph-fixed}, applied in the same way to
$E^*$-valued maps, identifies its fixed-point space with
$\BLip_0^{ph}(X,Y;E^*)$.  Formula
\eqref{eq:dual-valued-dilation-mean} is exactly
\eqref{eq:invariant-mean-projection} in this case.

It remains only to note that the range is nonzero when $X$, $Y$, and
$E$ are nonzero.  Choose nonzero $x^*\in X^*$, $y^*\in Y^*$, and
$e^*\in E^*$.  Then
$T(x,y)=x^*(x)y^*(y)e^*$
is a nonzero element of $\BLip_0^{ph}(X,Y;E^*)$.  Hence, the projection
has norm exactly one.
\end{proof}

\begin{corollary}\label{cor:norm-one-scalar-ph}
Let $X$ and $Y$ be nonzero Banach spaces.  Then there exists a norm-one
linear projection
$P_m:\BLip_0(X,Y)\longrightarrow\BLip_0^{ph}(X,Y)$.
More precisely,
\begin{equation}\label{eq:scalar-dilation-mean}
P_mT(x,y)
=
m_{(a,b)}\left[(a,b) \longmapsto
ab\,T\left(\frac{x}{a},\frac{y}{b}\right)
\right]
\end{equation}
and $\BLip(P_mT)\leq\BLip(T)$.
\end{corollary}

\begin{corollary}\label{cor:reflexive-codomain-ph}
If $F$ is reflexive and $X,Y$ are nonzero Banach spaces, then
$\BLip_0^{ph}(X,Y;F)$ is contractively complemented in
$\BLip_0(X,Y;F)$.  If this space is nonzero, the corresponding
projection has norm one.
\end{corollary}

\begin{proof}
Apply the bidual-valued averaging described in
Remark~\ref{rem:general-codomain-bidual} and use the canonical
surjective isometry $J_F:F\to F^{**}$.
\end{proof}

\begin{remark}\label{rem:two-projections}
Proposition~\ref{prop:radial-projection} and
Theorem~\ref{thm:amenable-projection} provide two complementary
projection mechanisms.  The radial projection is explicit, canonical,
and works for an arbitrary Banach-valued codomain, but the estimate
obtained above is \(\|P_{\mathrm{rad}}\|\leq9\).  If the codomain is a
dual space, amenability gives the invariant-mean projection \(P_m\) of
norm one.  The latter is nonconstructive and need not be $weak^*$
continuous.  Thus the radial construction reflects the orbit
cross-section \(\mathbb{S}_X\times \mathbb{S}_Y\), whereas invariant averaging reflects
amenability of the dilation group.

For a related one-variable norm-one projection from a Lipschitz space
onto its linear subspace in the dual-valued setting, see \cite{AM3}.
\end{remark}

\subsection{Concrete examples}

We finish with two examples illustrating that the size of
\(\BLip_0^{ph}(X,Y)\) depends strongly on the geometry of the domain
spaces.

\begin{proposition}\label{prop:R-R-four-dimensional}
The space
$\BLip_0^{ph}(\mathbb R,\mathbb R;\mathbb R)$
is four-dimensional, with basis
$\{T_{++},T_{+-},T_{-+},T_{--}\}$,
where, for \(z\in\mathbb R\),
$z^+=\max\{z,0\}, ~ \mbox{and}~z^-=\max\{-z,0\}$,
and  $$T_{++}(x,y)=x^+y^+, ~
T_{+-}(x,y)=x^+y^-,~T_{-+}(x,y)=x^-y^+,~\mbox{and}~
T_{--}(x,y)=x^-y^-.$$
\end{proposition}

\begin{proof}
Immediately, it follows that each of the four maps vanishes on the coordinate axes and is positively
homogeneous in each variable.  Moreover,
$$T_{++}(x,y)-T_{++}(x',y)-T_{++}(x,y')
+T_{++}(x',y')
=(x^+-x'^+)(y^+-y'^+),
$$
and therefore \(\BLip(T_{++})=1\), since \(t\mapsto t^+\) is
\(1\)-Lipschitz.  The same argument applies to the other three maps.

Let \(T\in\BLip_0^{ph}(\mathbb R,\mathbb R;\mathbb R)\).  Then from positive
homogeneity on the four sign quadrants we have
$$T(x,y)
=T(1,1)x^+y^+
+T(1,-1)x^+y^- + T(-1,1)x^-y^+
+T(-1,-1)x^-y^-.$$
Hence, the four displayed maps span the space.  Evaluating a linear
combination at
$(1,1),~(1,-1),$ $(-1,1)$ and $~(-1,-1)$,
shows that they are linearly independent.

In particular, the evaluation operator
$\Phi(T)
=
\bigl(
T(1,1),T(1,-1),T(-1,1),T(-1,-1)
\bigr)$
is a linear isomorphism onto \(\mathbb R^4\).  If
\(\Phi(T)=(a,b,c,d)\), then
$\|\Phi(T)\|_{\ell_\infty^4}
\leq
\BLip(T)
\leq
\|\Phi(T)\|_{\ell_1^4}
\leq
4\|\Phi(T)\|_{\ell_\infty^4}$.
Thus \(\Phi\) is a linear homeomorphism.
\end{proof}

\begin{remark}\label{rem:quadrants-Linfty}
Under the identification of
\(\BLip_0(\mathbb R,\mathbb R;\mathbb R)\) with
\(L^\infty(\mathbb R^2)\) used in \cite{AM1}, the four basis elements
above correspond to the indicator functions of the four open
quadrants.  Thus the preceding proposition identifies the positively
homogeneous subspace with the four-dimensional quadrant subspace of
\(L^\infty(\mathbb R^2)\).
\end{remark}

\begin{proposition}\label{prop:R2-R2-infinite-dimensional}
The space
$\BLip_0^{ph}(\mathbb R^2,\mathbb R^2;\mathbb R)$ is infinite-dimensional.
\end{proposition}

\begin{proof}
By Theorem~\ref{thm:sphere-restriction}, it is enough to construct an
infinite linearly independent family in
$\BLip_0(\mathbb{S}_{\mathbb R^2}^0,\mathbb{S}_{\mathbb R^2}^0;\mathbb R)$.
We identify \(\mathbb{S}_{\mathbb R^2}\) as
$\left\{(\cos\theta,\sin\theta): \theta \in [0,2\pi)\right\}$.  For each \(n\in\mathbb N\), we define $f_n:\mathbb{S}_{\mathbb R^2}^0\longrightarrow\mathbb R $ given by 
$$  
 f_n(u) = \begin{cases} \cos(n\theta), &\text{if }u=(\cos\theta,\sin\theta)\in \mathbb{S}_{\mathbb R^2},\\[1mm] 0, &\text{if }u=0. \end{cases} .$$
 
Therefore, each \(f_n\) is Lipschitz on \(\mathbb{S}_{\mathbb R^2}^0\).
Next, we define $g:\mathbb{S}_{\mathbb R^2}^0\longrightarrow\mathbb R$ given by $g(u)=u_1,~ u=(u_1,u_2)\in \mathbb{S}_{\mathbb R^2},$ and $g(0)=0$.
Using $f_n, g$ we further consider
$
G_n(u,v)=f_n(u)g(v)$. 
Then
$$
G_n(u,v)-G_n(u',v)-G_n(u,v')+G_n(u',v')=
\bigl(f_n(u)-f_n(u')\bigr)
\bigl(g(v)-g(v')\bigr),$$
so that $G_n \in \BLip_0(\mathbb{S}_{\mathbb R^2}^0,\mathbb{S}_{\mathbb R^2}^0;\mathbb R)$. 

We now show that the family $\{G_n:n\geq1\}$ is linearly independent. Suppose that, for some \(m\in\mathbb N\) and scalars \(c_1,\ldots,c_m\),~ $ \sum_{n=1}^m c_nG_n=0$. Now evaluating the summand at $ v=e_1=(1,0)\in \mathbb{S}_{\mathbb R^2}$, we obtain $ 0 = \sum_{n=1}^m c_nG_n(u,e_1) = \sum_{n=1}^m c_nf_n(u)$ for every \(u\in \mathbb{S}_{\mathbb R^2}\). Again, writing $ u=(\cos\theta,\sin\theta)$ we get $ \sum_{n=1}^m c_n\cos(n\theta)=0 $ for any $0\leq\theta<2\pi$. As, the functions $ \cos\theta,\cos(2\theta),\ldots,\cos(m\theta)$ are linearly independent, the family $\{G_n:n\geq1\}$ is linearly independent.

Therefore, applying the radial extension from
Theorem~\ref{thm:sphere-restriction} gives an infinite linearly
independent family in
\(\BLip_0^{ph}(\mathbb R^2,\mathbb R^2;\mathbb R)\).
\end{proof}

\section{Harmonic two-Lipschitz mappings on finitely generated groups}
\label{subsec:two-Lipschitz-harmonic}

The preceding application arises from the product dilation action. We now consider a different source of invariant pointwise-closed subspaces, arising from normalized translations on finitely generated groups. Although ordinary left translations do not preserve the distinguished identity, recentering induces natural isometric actions on the corresponding Lipschitz and two-Lipschitz spaces. Harmonicity in each variable is then expressed by pointwise linear constraints, and the general quotient-predual philosophy developed above yields a natural harmonic tensor predual.

\subsection{Finitely generated groups and word metrics}
\label{subsec:word-metric}

Let $G$ be a group with identity element $e_G$. A subset
$S\subseteq G$ is called a \emph{generating set} for $G$ if every
element $x\in G$ can be written as a finite product of elements of
$S\cup S^{-1}$. The group $G$ is said to be \emph{finitely generated}
if it admits a finite generating set. For the whole section we assume that $e_G \notin S$.

Throughout this section, all generating sets are assumed to be
\emph{symmetric}; that is,
$S=S^{-1}$.
Thus, if $S$ generates $G$, every $x\in G$ can be written in the form
$x=s_1s_2\cdots s_n,~ ~ \mbox{where}~ s_1,\ldots,s_n\in S$. For any $x\in G$, the \emph{word length} of $x$ with respect to $S$ is
defined by
$$|x|_S
=
\min\left\{
n\in\mathbb N:
x=s_1\cdots s_n,\;
s_1,\ldots,s_n\in S
\right\},$$
with the convention $|e_G|_S=0$.
The associated \emph{word metric} is
$d_S(x,y)
=
|x^{-1}y|_S,~\mbox{for any}~ x,y\in G$. Now, from the symmetry of $S$ we have
$|x^{-1}|_S=|x|_S$,
and hence $d_S$ is indeed a metric. We record below the elementary
properties that will be used later (the proof of it follows from simple computation).

\begin{proposition}\label{prop:word-metric-properties}
Let $G$ be a group generated by a finite symmetric set $S$. Then:

\begin{enumerate}
\item For all $x,y\in G~|xy|_S\leq |x|_S+|y|_S$.

\item The word metric is left invariant; namely,
$d_S(ax,ay)=d_S(x,y)$
for every $a,x,y\in G$.

\item For every $x\in G$ and $s\in S~d_S(x,xs)=d_S(x,xs^{-1})=1$.

\item If
$x^{-1}y=s_1s_2\cdots s_n,
~\mbox{where}~ s_1,\ldots,s_n\in S$,
is a representation of minimal length, then
$n=d_S(x,y)$,
and, setting
$x_0=x,\qquad
x_j=xs_1\cdots s_j
\quad(1\leq j\leq n)$,
we have $x_n=y~\text{and}~d_S(x_{j-1},x_j)=1$
for every $j=1,\ldots,n$.

\item Every closed ball
$B_S(e_G,n)
=
\{x\in G:|x|_S\leq n\}$
is finite.
\end{enumerate}
\end{proposition}

\begin{remark}\label{rem:right-multiplication-word-metric}
Although the word metric is always left invariant, right translations
need not be isometries in general. What will be important for us is
the elementary identity
$d_S(x,xs^{-1})=1,~
\mbox{for any}~ s\in S$,
which allows discrete right derivatives to detect the Lipschitz norm.
\end{remark}
\subsubsection{Lipschitz maps on finitely generated groups}

Let $G$ be a finitely generated group and let $S$ be a finite symmetric
generating set. Following \cite[Section~1.7.2, p.~31]{Yadin}, for a
function $f:G\to\mathbb R$ and $s\in S$, the \emph{right derivative of
$f$ in the direction $s$} is defined by
$\partial_s f(x)=f(xs^{-1})-f(x),~\mbox{for any}~ x\in G$.
The associated discrete gradient is
$\nabla_S f(x)
=
\bigl(\partial_s f(x)\bigr)_{s\in S}$,
and its supremum seminorm is
$\|\nabla_S f\|_\infty
=
\sup_{x\in G}\sup_{s\in S}
|\partial_s f(x)|$.
With respect to the word metric $d_S$, this seminorm coincides with the
usual Lipschitz seminorm; thus
$\Lip(f)=\|\nabla_S f\|_\infty$.
In particular, $f$ is Lipschitz if and only if
$\|\nabla_S f\|_\infty<\infty$. Moreover, the resulting class of
Lipschitz functions is independent of the choice of the finite symmetric
generating set; see \cite[Exercise~1.94, pp.~31--32]{Yadin}.

The preceding formulation provides the natural starting point for the
two-variable setting. Indeed, whereas the Lipschitz seminorm is detected
by first-order right differences, the defining inequality for a
two-Lipschitz mapping involves a mixed difference in the two variables.
This suggests replacing the discrete gradient above by a mixed discrete
gradient obtained by taking one right derivative in each variable.
\subsection{Mixed derivatives and the two-gradient}
Let $G$ and $H$ be finitely generated groups with finite symmetric
generating sets $S_G$ and $S_H$, respectively. We equip $G$ and $H$
with the corresponding word metrics $d_{S_{G}}, d_{S_{H}}$, respectively and denote their identity elements
by $e_G$ and $e_H$.

Following the usual discrete differential notation, for a function
$T:G\times H\longrightarrow\mathbb R$,
and $s\in G$, $t\in H$, we define the right derivatives in the first
and second variables by
$$\partial_s^{(1)}T(x,y)
=
T(xs^{-1},y)-T(x,y)~\mbox{and}~\partial_t^{(2)}T(x,y)
=
T(x,yt^{-1})-T(x,y).$$
Since these operators act on different variables, they commute. We
therefore define the mixed right derivative by
$\partial_{s,t}^{(1,1)}T
=
\partial_s^{(1)}\partial_t^{(2)}T$.
Explicitly,
$$\partial_{s,t}^{(1,1)}T(x,y)
=
T(xs^{-1},yt^{-1})
-T(xs^{-1},y)-T(x,yt^{-1})
+T(x,y).$$
\label{eq:mixed-right-derivative}
For any $(x,y)\in G\times H$, we define the \emph{two-gradient}
$\nabla^{(2)}T: G \times H \to \mathbb K^{S_G\times S_H}$ given by $\nabla^{(2)}T(x,y)
=
\bigl(
\partial_{s,t}^{(1,1)}T(x,y)
\bigr)_{(s,t)\in S_G\times S_H}$.
We put
\begin{equation}\label{eq:two-gradient-norm}
\|\nabla^{(2)}T\|_\infty
=
\sup_{\substack{x\in G,\ y\in H\\
s\in S_G,\ t\in S_H}}
\left|
\partial_{s,t}^{(1,1)}T(x,y)
\right|.
\end{equation}

\begin{proposition}\label{prop:gradient-BLip-norm}
Let $T: G \times H \to \mathbb{R}$ be such that
$T(e_G,y)=T(x,e_H)=0 $ for any $x\in G,\ y\in H$.
Then $T$ is a two-Lipschitz mapping if and only if
$\|\nabla^{(2)}T\|_\infty<\infty$.
Moreover,
$\BLip(T)
=
\|\nabla^{(2)}T\|_\infty$.
\end{proposition}

\begin{proof}
Suppose first that $T$ is two-Lipschitz. Since
$d_{S_{G}}(x,xs^{-1})=1,~\mbox{and}~d_{S_{H}}(y,yt^{-1})=1$.
For $s\in S_G$ and $t\in S_H$, the defining two-Lipschitz inequality
immediately gives
$|\partial_{s,t}^{(1,1)}T(x,y)|
\leq
\BLip(T)$.
Hence
$\|\nabla^{(2)}T\|_\infty
\leq
\BLip(T)$.

Conversely, we put
$C=\|\nabla^{(2)}T\|_\infty$.
Let $x,x'\in G$ and $y,y'\in H$. Since $S_G$ and $S_H$ are symmetric
generating sets, we may choose
$s_1,\ldots,s_m\in S_G
~\text{and}~
t_1,\ldots,t_n\in S_H$
such that
$x^{-1}x'=s_1^{-1}\cdots s_m^{-1},
~
y^{-1}y'=t_1^{-1}\cdots t_n^{-1},$
where
$m=d_{S_{G}}(x,x'),
~n=d_{S_{H}}(y,y')$.
We again consider
$x_0=x,~
x_i=x_{i-1}s_i^{-1} ~\mbox{for}~1\leq i\leq m,$
and $y_0=y,~y_j=y_{j-1}t_j^{-1}~ \mbox{for}~1\leq j\leq n$.
Then
$x_m=x',~
y_n=y'$.

We now expand the mixed difference by telescoping first in the
$y$-variable and then in the $x$-variable. For a fixed
$i\in\{1,\ldots,m\}$, we set
$A_j
=
T(x_{i-1},y_j)-T(x_i,y_j),
~\mbox{for}~
0\leq j\leq n$.
Then
\[
\begin{aligned}
\sum_{j=1}^n
&\Bigl[
T(x_{i-1},y_{j-1})
-T(x_i,y_{j-1})
-T(x_{i-1},y_j)
+T(x_i,y_j)
\Bigr]=
\sum_{j=1}^n(A_{j-1}-A_j)\\
&=
A_0-A_n=
T(x_{i-1},y)-T(x_i,y)-T(x_{i-1},y')
+T(x_i,y').
\end{aligned}
\]
Now, summing this identity over $i=1,\ldots,m$, we obtain
\[
\begin{aligned}
&\sum_{i=1}^m\sum_{j=1}^n
\Bigl[
T(x_{i-1},y_{j-1})
-T(x_i,y_{j-1})
-T(x_{i-1},y_j)
+T(x_i,y_j)
\Bigr]\\
&=
\sum_{i=1}^m
\Bigl[
T(x_{i-1},y)-T(x_i,y)
\Bigr]-
\sum_{i=1}^m
\Bigl[
T(x_{i-1},y')-T(x_i,y')
\Bigr].
\end{aligned}
\]
Both sums on the right telescope in the first variable. Hence
\[
\begin{aligned}
&\sum_{i=1}^m\sum_{j=1}^n
\Bigl[
T(x_{i-1},y_{j-1})
-T(x_i,y_{j-1})-T(x_{i-1},y_j)
+T(x_i,y_j)
\Bigr]=T(x,y)-T(x',y)-T(x,y')+T(x',y').
\end{aligned}
\]
Therefore,
\begin{equation}\label{eq:double-telescoping}
\begin{aligned}
&T(x,y)-T(x',y)-T(x,y')+T(x',y')=
\sum_{i=1}^m\sum_{j=1}^n
\Bigl[
T(x_{i-1},y_{j-1})
-T(x_i,y_{j-1})
-T(x_{i-1},y_j)
+T(x_i,y_j)
\Bigr].
\end{aligned}
\end{equation}

Since
$x_i=x_{i-1}s_i^{-1},~y_j=y_{j-1}t_j^{-1}$,
the summand in \eqref{eq:double-telescoping} is precisely
$\partial_{s_i,t_j}^{(1,1)}
T(x_{i-1},y_{j-1})$.
Hence, its absolute value is at most $C$. It follows that
\[
|T(x,y)-T(x',y)-T(x,y')+T(x',y')|\leq
\sum_{i=1}^m\sum_{j=1}^n
\left|
\partial_{s_i,t_j}^{(1,1)}
T(x_{i-1},y_{j-1})
\right|\leq
Cmn=
C\,d_G(x,x')\,d_H(y,y').
\]
Thus $T$ is two-Lipschitz and
$\BLip(T)\leq C
=
\|\nabla^{(2)}T\|_\infty$.
Combining this with the first inequality we have
$\BLip(T)
=
\|\nabla^{(2)}T\|_\infty$.
\end{proof}
\begin{proposition}[Independence of the generating sets]
\label{prop:two-gradient-generators}
Let $G$ and $H$ be finitely generated groups. Let
$S_G,S_G'$ be finite symmetric generating sets of $G$, and let
$S_H,S_H'$ be finite symmetric generating sets of $H$.

For a function
$T:G\times H\longrightarrow\mathbb R,$ we consider 
$\|\nabla^{(2)}_{S_G,S_H}T\|_\infty
=
\sup_{\substack{x\in G,\ y\in H\\
s\in S_G,\ t\in S_H}}
\left|
\partial_{s,t}^{(1,1)}T(x,y)
\right|$,
and analogously define
$\|\nabla^{(2)}_{S_G',S_H'}T\|_\infty$.
We put
$C_G
=
\max_{s\in S_G}|s|_{S_G'},~\mbox{and}~
C_H
=
\max_{t\in S_H}|t|_{S_H'}$.
Then
\begin{equation}\label{eq:gradient-generator-comparison}
\|\nabla^{(2)}_{S_G,S_H}T\|_\infty
\leq
C_GC_H\,
\|\nabla^{(2)}_{S_G',S_H'}T\|_\infty.
\end{equation}

Similarly, if
$C_G'
=
\max_{s'\in S_G'}|s'|_{S_G},~\mbox{and}~
C_H'
=
\max_{t'\in S_H'}|t'|_{S_H},$
then
\begin{equation}\label{eq:gradient-generator-comparison-reverse}
\|\nabla^{(2)}_{S_G',S_H'}T\|_\infty
\leq
C_G'C_H'\,
\|\nabla^{(2)}_{S_G,S_H}T\|_\infty.
\end{equation}
Consequently,
$\|\nabla^{(2)}_{S_G,S_H}T\|_\infty<\infty
\quad\Longleftrightarrow\quad
\|\nabla^{(2)}_{S_G',S_H'}T\|_\infty<\infty$.
Thus, the class of functions having bounded two-gradient is independent
of the choice of finite symmetric generating sets, and the corresponding
two-gradient seminorms are equivalent.
\end{proposition}

\begin{proof}
Since $S_G$ is finite and $S_G'$ generates $G$, every $s\in S_G$
admits a representation
$s=r_1r_2\cdots r_p,~\mbox{where}~
r_1,\ldots,r_p\in S_G'$
with
$p=|s|_{S_G'}\leq C_G$.
Likewise, for every $t\in S_H$, we may write
$t=q_1q_2\cdots q_\ell,~ \mbox{where}~q_1,\ldots,q_\ell\in S_H'$,
where $\ell=|t|_{S_H'}\leq C_H$.

Let $x\in G, ~y\in H,
~s\in S_G, ~t\in S_H$.
Since
$s^{-1}=r_p^{-1}\cdots r_1^{-1},~
t^{-1}=q_\ell^{-1}\cdots q_1^{-1}$,
we define
$x_0=x,~x_i=x_{i-1}r_{p-i+1}^{-1},~\mbox{for}~ 1\leq i\leq p$,
and
$y_0=y,~y_j=y_{j-1}q_{\ell-j+1}^{-1},~\mbox{for}~ 1\leq j\leq \ell$.
Then
$x_p=xs^{-1},~\mbox{and}~
y_\ell=yt^{-1}$.

Again, by the double telescoping identity,
\begin{align*}
&T(xs^{-1},yt^{-1})
-T(xs^{-1},y)
-T(x,yt^{-1})
+T(x,y)
=
\sum_{i=1}^{p}\sum_{j=1}^{\ell}
\Bigl[
T(x_i,y_j)-T(x_i,y_{j-1})
-T(x_{i-1},y_j)+T(x_{i-1},y_{j-1})
\Bigr].
\end{align*}
For each $i$ and $j$, we put
$a_i=r_{p-i+1}\in S_G',~\mbox{and}~
b_j=q_{\ell-j+1}\in S_H'$.
Since $x_i=x_{i-1}a_i^{-1},~\mbox{and}~
y_j=y_{j-1}b_j^{-1}$,
the corresponding summand is exactly
$\partial_{a_i,b_j}^{(1,1)}
T(x_{i-1},y_{j-1})$.
Therefore,
\begin{align*}
\left|
\partial_{s,t}^{(1,1)}T(x,y)
\right|
\leq
\sum_{i=1}^{p}\sum_{j=1}^{\ell}
\left|
\partial_{a_i,b_j}^{(1,1)}
T(x_{i-1},y_{j-1})
\right|
\leq
p\ell\,
\|\nabla^{(2)}_{S_G',S_H'}T\|_\infty
\leq
C_GC_H\,
\|\nabla^{(2)}_{S_G',S_H'}T\|_\infty.
\end{align*}
Now, taking the supremum over
$x\in G,~ y\in H,~
s\in S_G,~ t\in S_H$ we have
$\|\nabla^{(2)}_{S_G,S_H}T\|_\infty
\leq
C_GC_H\,
\|\nabla^{(2)}_{S_G',S_H'}T\|_\infty$.

Interchanging $(S_G,S_H)$ with $(S_G',S_H')$ gives
$\|\nabla^{(2)}_{S_G',S_H'}T\|_\infty
\leq
C_G'C_H'\,
\|\nabla^{(2)}_{S_G,S_H}T\|_\infty$.
Hence, the two seminorms are equivalent, and in particular one is finite
if and only if the other is finite.
\end{proof}

\begin{lemma}\label{lem:zero-right-derivative}
Let $G$ be a group generated by a set $S$, and let
$f:G\longrightarrow\mathbb R$ be a function. Suppose that
$\partial_s f(x)
=
f(xs^{-1})-f(x)
=
0$
for every $x\in G$ and every $s\in S$. Then $f$ is constant on $G$.
\end{lemma}

\begin{proof}
Let $x\in G$. Since $S$ generates $G$, there exist
$s_1,\ldots,s_n\in S\cup S^{-1}$ such that
$x=s_1s_2\cdots s_n$.
The assumption
$f(zs^{-1})=f(z)~\mbox{for any}~z\in G,\ s\in S$
also implies $f(zs)=f(z)$.
Hence, $f$ is unchanged under right multiplication by any element of
$S\cup S^{-1}$. Applying this successively along the above word gives
$f(x)
=
f(s_1\cdots s_n)
=
f(s_1\cdots s_{n-1})
=
\cdots
=
f(e_G)$.
Thus, $f$ is constant.
\end{proof}

\begin{proposition}\label{prop:zero-two-gradient}
Let $G$ and $H$ be groups generated by $S_G$ and $S_H$, respectively,
and let
$T:G\times H\longrightarrow\mathbb R$
be a function. The following assertions are equivalent:

\begin{enumerate}
\item
$\partial_{s,t}^{(1,1)}T(x,y)=0$
for every $x\in G$, $y\in H$, $s\in S_G$, and $t\in S_H$; equivalently,
$\|\nabla^{(2)}T\|_\infty=0$.

\item There exist functions
$f:G\longrightarrow\mathbb R,~
g:H\longrightarrow\mathbb R$
such that
$T(x,y)=f(x)+g(y)~\mbox{for any}~x\in G,\ y\in H$.
\end{enumerate}

More precisely, whenever these conditions hold,
\begin{equation}\label{eq:zero-mixed-gradient-decomposition}
T(x,y)
=
T(x,e_H)+T(e_G,y)-T(e_G,e_H)
\end{equation}
for every $x\in G$ and $y\in H$.
Consequently, on the normalized class of functions satisfying
$T(e_G,y)=0~\text{and}~
T(x,e_H)=0$
for every $x\in G$ and $y\in H$, the seminorm
$T\longmapsto\|\nabla^{(2)}T\|_\infty$
is a norm.
\end{proposition}

\begin{proof}
Let us assume first that
$T(x,y)=f(x)+g(y)$.
For $s\in S_G$ and $t\in S_H$, we have
\begin{align*}
\partial_{s,t}^{(1,1)}T(x,y)
={}&
T(xs^{-1},yt^{-1})
-T(xs^{-1},y)-T(x,yt^{-1})
+T(x,y)\\
={}&
f(xs^{-1})+g(yt^{-1})
-f(xs^{-1})-g(y)-f(x)-g(yt^{-1})
+f(x)+g(y)\\
={}&0.
\end{align*}
Hence,
$\|\nabla^{(2)}T\|_\infty=0$.

Conversely, suppose that $\partial_{s,t}^{(1,1)}T(x,y)=0$
for all $x\in G$, $y\in H$, $s\in S_G$, and $t\in S_H$.
Let, $t\in S_H$ and $y\in H$. We define
$\varphi_{t,y}:G\longrightarrow\mathbb R$
by
$\varphi_{t,y}(x)
=
\partial_t^{(2)}T(x,y)
=
T(x,yt^{-1})-T(x,y)$.
Now, for any $s\in S_G$,
\begin{align*}
\partial_s\varphi_{t,y}(x)
=
\varphi_{t,y}(xs^{-1})-\varphi_{t,y}(x)=
T(xs^{-1},yt^{-1})
-T(xs^{-1},y)-T(x,yt^{-1})
+T(x,y)=
\partial_{s,t}^{(1,1)}T(x,y)=0.
\end{align*}
By Lemma~\ref{lem:zero-right-derivative}, $\varphi_{t,y}$ is constant
on $G$. Therefore,
\begin{equation}\label{eq:mixed-gradient-independent-x}
\partial_t^{(2)}T(x,y)
=
\partial_t^{(2)}T(e_G,y)
\end{equation}
for every $x\in G$, $y\in H$, and $t\in S_H$.

Now for fixed $x\in G$ and we define
$F_x:H\longrightarrow\mathbb R$
by
$F_x(y)
=
T(x,y)-T(e_G,y)$.
Using \eqref{eq:mixed-gradient-independent-x}, for every $t\in S_H$, we have
\begin{align*}
\partial_tF_x(y)
&=
F_x(yt^{-1})-F_x(y)\\
&=
T(x,yt^{-1})-T(x,y)-\bigl(T(e_G,yt^{-1})-T(e_G,y)\bigr)\\
&=
\partial_t^{(2)}T(x,y)
-
\partial_t^{(2)}T(e_G,y)\\
&=0.
\end{align*}
Again by Lemma~\ref{lem:zero-right-derivative}, now applied to $H$,
the function $F_x$ is constant on $H$. Consequently,
$F_x(y)=F_x(e_H)$
for every $y\in H$. Thus
$T(x,y)-T(e_G,y)
=
T(x,e_H)-T(e_G,e_H)$,
which is equivalent to
$T(x,y)
=
T(x,e_H)+T(e_G,y)-T(e_G,e_H)$.
This proves \eqref{eq:zero-mixed-gradient-decomposition}.

Defining
$f(x)=T(x,e_H)-T(e_G,e_H)$
and
$g(y)=T(e_G,y)$,
we obtain $T(x,y)=f(x)+g(y)$,
and hence, $(1)$ implies $(2)$.

Finally, suppose that $T$ belongs to the normalized class and $\|\nabla^{(2)}T\|_\infty=0$.
By \eqref{eq:zero-mixed-gradient-decomposition},
\[
T(x,y)
=
T(x,e_H)+T(e_G,y)-T(e_G,e_H)
=
0
\]
for every $(x,y)\in G\times H$. Thus, the only normalized function with
zero two-gradient seminorm is the zero function. Hence
$\|\nabla^{(2)}\cdot\|_\infty$ is a norm on the normalized class.
\end{proof}

The discrete differential formulation above also provides a natural
motivation for harmonicity. Given finitely supported probability measures
$\mu$ on $G$ and $\nu$ on $H$, one may average the right derivatives in
the corresponding variables. Thus, instead of requiring each discrete
directional derivative to vanish separately, one may require its
$\mu$-weighted or $\nu$-weighted average to vanish. Since the mixed
derivatives satisfy
\[
\partial_{s,t}^{(1,1)}
=
\partial_s^{(1)}\partial_t^{(2)},
\]
the two-gradient records the interaction between the discrete differential
structures in the two variables. This leads naturally to the notion of
two-Lipschitz mappings which are harmonic separately in each variable.
\subsection{Harmonicity in each variable} \label{harmonicity}

Let $\mu\in\ell_1(G)$ and $\nu\in\ell_1(H)$ be finitely supported
probability measures on $G$ and $H$, respectively. For a function $T:G\times H\to\mathbb \mathbb{R}$, we define
the following operators in the two variables by
$(P_\mu^{(1)}T)(x,y)
=
\sum_{s\in G}\mu(s)T(xs^{-1},y)$
and
$(P_\nu^{(2)}T)(x,y)
=
\sum_{t\in H}\nu(t)T(x,yt^{-1})$.
Equivalently, we put
$\Delta_\mu^{(1)}
=
P_\mu^{(1)}-I,
~\mbox{and}~
\Delta_\nu^{(2)}
=
P_\nu^{(2)}-I$.

In terms of the right derivatives introduced above, these operators satisfy
\[
\Delta_\mu^{(1)}T=(P_\mu^{(1)}-I)T
=
\sum_{s\in G}\mu(s)\,\partial_s^{(1)}T,
\qquad
\Delta_\nu^{(2)}T=(P_\nu^{(2)}-I)T
=
\sum_{t\in H}\nu(t)\,\partial_t^{(2)}T.
\]
Hence harmonicity in each variable is equivalent to the vanishing of the
corresponding weighted average of discrete right derivatives.
\begin{definition}\label{def:harmonic-two-Lipschitz}
A map
$T\in\BLip_0(G,H)$
is called \emph{$(\mu,\nu)$-harmonic in each variable} if
$\Delta_\mu^{(1)}T=0
~\text{and}~
\Delta_\nu^{(2)}T=0$.
Equivalently,
\begin{equation}\label{eq:mu-harmonic-first}
T(x,y)
=
\sum_{s\in G}\mu(s)T(xs^{-1},y)
\end{equation}
and
\begin{equation}\label{eq:nu-harmonic-second}
T(x,y)
=
\sum_{t\in H}\nu(t)T(x,yt^{-1})
\end{equation}
for every $x\in G$ and $y\in H$.
We denote the corresponding space by
$\HBLipz^{\mu,\nu}(G,H)$.
\end{definition}
Thus,
$\HBLipz^{\mu,\nu}(G,H)
=
\ker\Delta_\mu^{(1)}
\cap
\ker\Delta_\nu^{(2)}
\cap
\BLip_0(G,H)$.  

\begin{remark}[A spectral viewpoint]\label{rem:spectral-viewpoint}
The preceding harmonic construction is related to the
character-equivariant framework through a spectral viewpoint, although
the corresponding averaging operators arise from different actions.

Suppose that groups $G$ and $H$ act on pointed metric spaces $M$ and
$N$, respectively, by similarities, with associated characters
$\chi_G$ and $\chi_H$. Let $\mu$ and $\nu$ be finitely supported
probability measures on $G$ and $H$, respectively. We define
$(A_\mu^{(1)}T)(x,y)
=
\sum_{g\in G}\mu(g)T(gx,y),
~~
(A_\nu^{(2)}T)(x,y)
=
\sum_{h\in H}\nu(h)T(x,hy)$.
For $T\in \BLipz^{G\times H,\chi}(M,N),$ we have then $T(gx,hy)
=
\chi_G(g)\chi_H(h)T(x,y)$,
and therefore
$A_\mu^{(1)}T
=
\widehat\mu(\chi_G)\,T,
~~
A_\nu^{(2)}T
=
\widehat\nu(\chi_H)\,T$,
where
\[
\widehat\mu(\chi_G)
=
\sum_{g\in G}\mu(g)\chi_G(g),
~~
\widehat\nu(\chi_H)
=
\sum_{h\in H}\nu(h)\chi_H(h).
\]
Thus character-equivariant two-Lipschitz mappings are simultaneous
eigenvectors of the averaging operators $A_\mu^{(1)}$ and
$A_\nu^{(2)}$. In particular, if
$\widehat\mu(\chi_G)
=
\widehat\nu(\chi_H)
=
1$,
then $T$ is fixed by both averaging operators.

In the harmonic setting considered above, the operators
$P_\mu^{(1)}$ and $P_\nu^{(2)}$ are instead defined by the averages
$(P_\mu^{(1)}T)(x,y)
=
\sum_{s\in G}\mu(s)T(xs^{-1},y),
~~
(P_\nu^{(2)}T)(x,y)
=
\sum_{t\in H}\nu(t)T(x,yt^{-1})$.
Hence
$\HBLipz^{\mu,\nu}(G,H)
=
\left\{
T\in\BLip_0(G,H):
P_\mu^{(1)}T=T,\;
P_\nu^{(2)}T=T
\right\}$.
Thus both the character-equivariant and harmonic theories admit a
common spectral interpretation in terms of simultaneous eigenvectors
of two commuting averaging operators, while the underlying averaging
actions are different.
\end{remark}

\begin{proposition}\label{prop:commuting-harmonic-operators}
The operators $P_\mu^{(1)}$ and $P_\nu^{(2)}$ commute. More precisely,
for every function $T:G\times H\to\mathbb R$,
$P_\mu^{(1)}P_\nu^{(2)}T
=
P_\nu^{(2)}P_\mu^{(1)}T
=
P_{\mu\otimes\nu}T,$
where
$(P_{\mu\otimes\nu}T)(x,y)
=
\sum_{s\in G}\sum_{t\in H}
\mu(s)\nu(t)
T(xs^{-1},yt^{-1})$.
\end{proposition}

\begin{proof}
Let $(x,y)\in G\times H$. Since $\mu$ and $\nu$ are finitely
supported,
\begin{align*}
(P_\mu^{(1)}P_\nu^{(2)}T)(x,y)
=
\sum_{s\in G}
\mu(s)
(P_\nu^{(2)}T)(xs^{-1},y)=
\sum_{s\in G}\sum_{t\in H}
\mu(s)\nu(t)
T(xs^{-1},yt^{-1}).
\end{align*}
Similarly,
\begin{align*}
(P_\nu^{(2)}P_\mu^{(1)}T)(x,y)
=
\sum_{t\in H}
\nu(t)
(P_\mu^{(1)}T)(x,yt^{-1})=
\sum_{t\in H}\sum_{s\in G}
\nu(t)\mu(s)
T(xs^{-1},yt^{-1}).
\end{align*}
Since both sums are finite, their order may be interchanged. Hence
$P_\mu^{(1)}P_\nu^{(2)}
=
P_\nu^{(2)}P_\mu^{(1)}
=
P_{\mu\otimes\nu}$.

Thus, in the ambient function space,
\[
\HBLipz^{\mu,\nu}(G,H)
=
\BLip_0(G,H)\cap
E_{(1,1)}(P_\mu^{(1)},P_\nu^{(2)}),
\]
where $E_{(\cdot,\cdot)}(\cdot,\cdot)$ is the joint eigenspace.
\end{proof}
For preliminaries on joint eigenspace we refer \cite{JLT}.
\begin{corollary}\label{cor:separate-product-harmonic}
Every $(\mu,\nu)$-harmonic two-Lipschitz map is
$(\mu\otimes\nu)$-harmonic on the product group $G\times H$.
More precisely, if
$T\in \HBLipz^{\mu,\nu}(G,H)$,
then
$T(x,y)
=
\sum_{s\in G}\sum_{t\in H}
\mu(s)\nu(t)
T(xs^{-1},yt^{-1})$
for every $(x,y)\in G\times H$.
\end{corollary}

\begin{proof}
Since
$P_\mu^{(1)}T=T
~\text{and}~
P_\nu^{(2)}T=T$. Then
Proposition~\ref{prop:commuting-harmonic-operators} gives
$P_{\mu\otimes\nu}T
=
P_\mu^{(1)}P_\nu^{(2)}T
=
P_\mu^{(1)}T
=
T$.
\end{proof}

\begin{remark}\label{rem:product-vs-separate-harmonic}
Separate harmonicity is stronger than harmonicity with respect to the
product measure. Indeed,
$P_\mu^{(1)}T=T
\quad\text{and}\quad
P_\nu^{(2)}T=T$
imply
$P_{\mu\otimes\nu}T=T$,
whereas the equation
$P_\mu^{(1)}P_\nu^{(2)}T=T$
does not, in general, imply that $T$ is fixed separately by
$P_\mu^{(1)}$ and $P_\nu^{(2)}$.
Accordingly, throughout this section, by a harmonic two-Lipschitz map
we mean a map which is harmonic in each variable.
\end{remark}

\subsubsection{The normalized translation action}

As for $(a,b)\in G\times H$, ordinary left translation does not preserve
the conditions
$T(e_G,y)=T(x,e_H)=0$.
We therefore use its normalized version and define
\begin{align}
(U_{a,b}T)(x,y)
=
T(a^{-1}x,b^{-1}y)
-T(a^{-1}x,b^{-1})-T(a^{-1},b^{-1}y)
+T(a^{-1},b^{-1}),
\label{eq:recentered-translation}
\end{align}
Then
$(U_{a,b}T)(e_G,y)
=
(U_{a,b}T)(x,e_H)
=
0$.

\begin{proposition}\label{prop:translation-two-gradient}
For every $(a,b)\in G\times H$ the map $U_{a,b}:
\BLip_0(G,H)\longrightarrow\BLip_0(G,H)$ is a surjective linear isometry. Moreover,
$\|\nabla^{(2)}U_{a,b}T\|_\infty
=
\|\nabla^{(2)}T\|_\infty$, and the family $U:=(U_{a,b})_{(a,b)\in G\times H}$ defines an action of
$G\times H$ on $\BLip_0(G,H)$.
\end{proposition}

\begin{proof}
We first prove the invariance of the two-gradient norm. Since the mixed derivative annihilates every function depending on only one variable, the three normalizing terms in the definition of $U_{a,b}T$ do not contribute to the mixed derivative. Hence,
 a simple computation gives $\partial_{s,t}^{(1,1)}(U_{a,b}T)(x,y)
=
\partial_{s,t}^{(1,1)}T(a^{-1}x,b^{-1}y)$.
By taking suprema we have
$\|\nabla^{(2)}U_{a,b}T\|_\infty
=
\|\nabla^{(2)}T\|_\infty$.
Then, by Proposition~\ref{prop:gradient-BLip-norm},
$\BLip(U_{a,b}T)=\BLip(T)$.
From a direct computation it follows
$U_{a,b}U_{c,d}
=
U_{ac,bd}$,
and
$U_{a,b}^{-1}
=
U_{a^{-1},b^{-1}}$.
Thus, the family defines an action by surjective linear isometries.
\end{proof}

The preceding proposition is the two-variable analogue of the
translation invariance of the Lipschitz seminorm.

\begin{proposition}\label{prop:HBLip-invariant-pointwise}
Let $G$ and $H$ be finitely generated groups, and let $\mu$ and $\nu$
be finitely supported probability measures on $G$ and $H$, respectively.
Then $\HBLipz^{\mu,\nu}(G,H)$
is a linear, pointwise-closed subspace of $\BLip_0(G,H)$ and is invariant
under the normalized translation action
$(U_{a,b})_{(a,b)\in G\times H}$.
More precisely,
$U_{a,b}\bigl(\HBLipz^{\mu,\nu}(G,H)\bigr)
=
\HBLipz^{\mu,\nu}(G,H)$
for every $(a,b)\in G\times H$.

Consequently,
$\HBLipz^{\mu,\nu}(G,H)$
is weak$^*$ closed in
$\BLip_0(G,H)
=
\bigl(\F(G)\widehat\otimes_\pi\F(H)\bigr)^*$.
\end{proposition}
\begin{proof}
We first recall that
$\HBLipz^{\mu,\nu}(G,H)
=
\left\{
T\in\BLip_0(G,H):
P_\mu^{(1)}T=T
\ \text{and}\
P_\nu^{(2)}T=T
\right\}$.
Since both $P_\mu^{(1)}$ and $P_\nu^{(2)}$ are linear, it follows
immediately that $\HBLipz^{\mu,\nu}(G,H)$ is a linear subspace of
$\BLip_0(G,H)$.

We first prove that this subspace is closed for pointwise convergence.
Let $(T_\alpha)_\alpha$ be a net in
$\HBLipz^{\mu,\nu}(G,H)$ converging pointwise to some
$T\in\BLip_0(G,H)$. Thus,
$T_\alpha(x,y)\longrightarrow T(x,y)$
for every $(x,y)\in G\times H$.

Since each $T_\alpha$ is $\mu$-harmonic in the first variable, we have
$T_\alpha(x,y)
=
\sum_{s\in G}
\mu(s)T_\alpha(xs^{-1},y)$.
Since, $\mu$ is finitely supported, we may pass to the pointwise limit
inside the finite sum and obtain
\begin{align*}
T(x,y)=
\lim_\alpha T_\alpha(x,y)=
\lim_\alpha
\sum_{s\in G}
\mu(s)T_\alpha(xs^{-1},y)=
\sum_{s\in G}
\mu(s)
\lim_\alpha T_\alpha(xs^{-1},y)=
\sum_{s\in G}
\mu(s)T(xs^{-1},y).
\end{align*}
Hence
$P_\mu^{(1)}T=T$.
Similarly, since each $T_\alpha$ is $\nu$-harmonic in the second
variable and $\nu$ is finitely supported,
$T(x,y)
=
\sum_{t\in H}
\nu(t)T(x,yt^{-1})$,
and therefore
$P_\nu^{(2)}T=T$.
Thus
$T\in\HBLipz^{\mu,\nu}(G,H)$,
showing that $\HBLipz^{\mu,\nu}(G,H)$ is pointwise closed.

We next prove invariance under the normalized translation action. We again recall
that, 
\begin{align}
(U_{a,b}T)(x,y)
=T(a^{-1}x,b^{-1}y)
-T(a^{-1}x,b^{-1}) -T(a^{-1},b^{-1}y)
+T(a^{-1},b^{-1}).
\label{eq:normalized-translation-harmonic}
\end{align}
for any $(a,b)\in G\times H$.
Let
$T\in\HBLipz^{\mu,\nu}(G,H)$
and fix $(a,b)\in G\times H$. We first verify $\mu$-harmonicity of
$U_{a,b}T$ in the first variable. For $x\in G$ and $y\in H$,
\begin{align*}
(P_\mu^{(1)}U_{a,b}T)(x,y)=
\sum_{s\in G}
\mu(s)(U_{a,b}T)(xs^{-1},y)=
&\sum_{s\in G}\mu(s)
T(a^{-1}xs^{-1},b^{-1}y)
-
\sum_{s\in G}\mu(s)
T(a^{-1}xs^{-1},b^{-1})\\
&\qquad
-
\sum_{s\in G}\mu(s)
T(a^{-1},b^{-1}y)
+
\sum_{s\in G}\mu(s)
T(a^{-1},b^{-1}).
\end{align*}
Since $\mu$ is a probability measure,
$\sum_{s\in G}\mu(s)=1$.
Moreover, since $T$ is $\mu$-harmonic in the first variable,
$$\sum_{s\in G}
\mu(s)T(a^{-1}xs^{-1},b^{-1}y)
=
T(a^{-1}x,b^{-1}y)
~~\mbox{and} ~~\sum_{s\in G}
\mu(s)T(a^{-1}xs^{-1},b^{-1})
=
T(a^{-1}x,b^{-1}).$$
Therefore,
\begin{align*}
(P_\mu^{(1)}U_{a,b}T)(x,y)
=
T(a^{-1}x,b^{-1}y)
-T(a^{-1}x,b^{-1})-T(a^{-1},b^{-1}y)
+T(a^{-1},b^{-1})=(U_{a,b}T)(x,y).
\end{align*}
Thus
$P_\mu^{(1)}(U_{a,b}T)=U_{a,b}T$.
Hence $U_{a,b}T$ is $\mu$-harmonic in the first variable.

We now verify harmonicity in the second variable. For $x\in G$ and
$y\in H$,
\begin{align*}
(P_\nu^{(2)}U_{a,b}T)(x,y)=
\sum_{t\in H}
\nu(t)(U_{a,b}T)(x,yt^{-1})&=
\sum_{t\in H}\nu(t)
T(a^{-1}x,b^{-1}yt^{-1})
-
\sum_{t\in H}\nu(t)
T(a^{-1}x,b^{-1})\\
&-
\sum_{t\in H}\nu(t)
T(a^{-1},b^{-1}yt^{-1})
+
\sum_{t\in H}\nu(t)
T(a^{-1},b^{-1}).
\end{align*}
Since $\nu$ is a probability measure,
$\sum_{t\in H}\nu(t)=1$.
Furthermore, $\nu$-harmonicity of $T$ in the second variable gives
$$\sum_{t\in H}
\nu(t)T(a^{-1}x,b^{-1}yt^{-1})
=
T(a^{-1}x,b^{-1}y)~~\mbox{
and}
~~\sum_{t\in H}
\nu(t)T(a^{-1},b^{-1}yt^{-1})
=
T(a^{-1},b^{-1}y).$$
Consequently,
\begin{align*}
(P_\nu^{(2)}U_{a,b}T)(x,y)
=
T(a^{-1}x,b^{-1}y)
-T(a^{-1}x,b^{-1})-T(a^{-1},b^{-1}y)
+T(a^{-1},b^{-1})
=
(U_{a,b}T)(x,y).
\end{align*}
Hence
$P_\nu^{(2)}(U_{a,b}T)=U_{a,b}T$.
Therefore, $U_{a,b}T$ is $\nu$-harmonic in the second variable.

It follows that
$U_{a,b}T\in\HBLipz^{\mu,\nu}(G,H)$,
and hence
$U_{a,b}
\bigl(\HBLipz^{\mu,\nu}(G,H)\bigr)
\subseteq
\HBLipz^{\mu,\nu}(G,H)$.
Since $(U_{a,b})_{(a,b)\in G\times H}$ is a group action and
$U_{a,b}^{-1}=U_{a^{-1},b^{-1}}$,
applying the preceding inclusion to $(a^{-1},b^{-1})$ gives the reverse
inclusion. Therefore
\[
U_{a,b}
\bigl(\HBLipz^{\mu,\nu}(G,H)\bigr)
=
\HBLipz^{\mu,\nu}(G,H).
\]

Thus $\HBLipz^{\mu,\nu}(G,H)$ is an invariant pointwise-closed linear
subspace of $\BLip_0(G,H)$. By
Theorem~\ref{thm:pointwise-invariant-BLip},
it is weak$^*$ closed in
$\BLip_0(G,H)
=
\bigl(\F(G)\widehat\otimes_\pi\F(H)\bigr)^*$.
\end{proof}

\subsubsection{The harmonic tensor predual} \label{har ten prdu}

We define the one-variable harmonic relation spaces by
\[
N_\mu(G)
=
\overline{\operatorname{span}}
\left\{
\delta_x-
\sum_{s\in G}\mu(s)\delta_{xs^{-1}}
:
x\in G
\right\}
\subseteq\F(G)~
\mbox{and}~
N_\nu(H)
=
\overline{\operatorname{span}}
\left\{
\delta_y-
\sum_{t\in H}\nu(t)\delta_{yt^{-1}}
:
y\in H
\right\}
\subseteq\F(H).
\]
We consider $\F_\mu^h(G)=\F(G)/N_\mu(G),
~~\F_\nu^h(H)=\F(H)/N_\nu(H)$.

The following result may be deduced from \cite{DouchaIPC}. We nevertheless include the proof for completeness, as it provides a natural motivation for the two-variable framework considered in the sequel.
\begin{proposition}\label{prop:one-variable-harmonic-predual}
Let $G$ be a finitely generated group equipped with a word metric, and
let $\mu$ be a finitely supported probability measure on $G$. We consider
$\LHFz^{\mu}(G)
=
\left\{
f\in\Lip_0(G):
P_\mu f=f
\right\}$; the space of \emph{harmonic Lipschitz maps} in one variable,
where
$(P_\mu f)(x)
=
\sum_{s\in G}\mu(s)f(xs^{-1})$.
Then
$\LHFz^{\mu}(G)
=
N_\mu(G)^\perp$
under the canonical identification
$\Lipz(G)=\F(G)^*$.
Consequently, there is a canonical isometric identification
\begin{equation}\label{eq:one-variable-harmonic-predual}
\LHFz^{\mu}(G)
\cong
\bigl(\F_\mu^h(G)\bigr)^*.
\end{equation}
\end{proposition}

\begin{proof}
We recall that every $f\in\Lipz(G)$ is identified with the bounded
linear functional $\widehat{f}$ on $\F(G)$ determined by
$\widehat{f}(\delta_x)=f(x),
~\mbox{for any}~ x\in G$ (from the universal property).
Therefore, for every $x\in G$,
\begin{align*}
\widehat{
f}\left(
\delta_x-\sum_{s\in G}\mu(s)\delta_{xs^{-1}}
\right)
=f(x)-\sum_{s\in G}\mu(s)f(xs^{-1})=
f(x)-(P_\mu f)(x).
\end{align*}
Hence, $f$ annihilates the generating relations of $N_\mu(G)$ if and
only if
$P_\mu f=f$.
By continuity, this is equivalent to
$f\in N_\mu(G)^\perp$.
Thus
$\LHFz^{\mu}(G)=N_\mu(G)^\perp$.
Since $N_\mu(G)$ is norm closed, from the standard duality for quotient
spaces we have
$\bigl(\F(G)/N_\mu(G)\bigr)^*
\cong
N_\mu(G)^\perp$. Therefore, $\LHFz^{\mu}(G)
\cong
\bigl(\F_\mu^h(G)\bigr)^*$.
\end{proof}
The same argument gives
$\LHFz^{\nu}(H)
\cong
\bigl(\F_\nu^h(H)\bigr)^*$.

We now define
\begin{equation}\label{eq:harmonic-tensor-relation}
K_{\mu,\nu}
=
\overline{
N_\mu(G)\otimes\F(H)
+
\F(G)\otimes N_\nu(H)
}^{\,\|\cdot\|_\pi}.
\end{equation}

\begin{proposition}\label{prop:harmonic-tensor-predual}
Let $G$ and $H$ be finitely generated groups, and let
$\mu$ and $\nu$ be finitely supported probability measures on $G$
and $H$, respectively. 
Then
$\HBLipz^{\mu,\nu}(G,H)
=
K_{\mu,\nu}^{\perp}$.
Consequently, there is a canonical isometric identification
\begin{equation}\label{eq:HBLip-quotient-predual}
\HBLipz^{\mu,\nu}(G,H)
\cong
\left(
\bigslant{
\F(G)\widehat{\otimes}_\pi\F(H)
}{
K_{\mu,\nu}
}
\right)^*.
\end{equation}

Moreover,
\begin{equation}\label{eq:harmonic-tensor-factorization}
\HBLipz^{\mu,\nu}(G,H)
\cong
\left(
\F_\mu^h(G)
\widehat{\otimes}_\pi
\F_\nu^h(H)
\right)^*
\end{equation}
canonically and isometrically.
\end{proposition}

\begin{proof}
Let
$T\in\BLip_0(G,H)$
and $T_L:
\F(G)\widehat{\otimes}_\pi\F(H)
\longrightarrow \mathbb R$
denote its canonical linearization.

We first consider harmonicity in the first variable. For
$x\in G$ and $y\in H$, the identity
$T(x,y)
=
\sum_{s\in G}\mu(s)T(xs^{-1},y)
$ is equivalent to
$T_L\left(
\left(
\delta_x-\sum_{s\in G}\mu(s)\delta_{xs^{-1}}
\right)\otimes\delta_y
\right)
=0$.
Hence $T$ is $\mu$-harmonic in the first variable if and only if
$T_L$ vanishes on
$N_\mu(G)\otimes\F(H)$.

Similarly, $T$ is $\nu$-harmonic in the second variable if and only if
$T_L$ vanishes on
$\F(G)\otimes N_\nu(H)$.
Therefore, using continuity of $T_L$ we have $T$ is harmonic in each variable if and only if
$T_L\big|_{K_{\mu,\nu}}=0$.
Moreover, under the canonical identification
we get $\BLip_0(G,H)
\cong
\left(
\F(G)\widehat{\otimes}_\pi\F(H)
\right)^*$,
this is precisely the condition
$T\in K_{\mu,\nu}^{\perp}$.
Thus
$\HBLipz^{\mu,\nu}(G,H)
=
K_{\mu,\nu}^{\perp}$.

The standard duality between quotients and annihilators now yields
$K_{\mu,\nu}^{\perp}
\cong
\left(
\bigslant{
\F(G)\widehat{\otimes}_\pi\F(H)
}{
K_{\mu,\nu}
}
\right)^*$,
canonically and isometrically, which proves
\eqref{eq:HBLip-quotient-predual}.

Finally, let
$q_\mu:\F(G)\longrightarrow\F_\mu^h(G),
~~
q_\nu:\F(H)\longrightarrow\F_\nu^h(H)$
be the canonical quotient maps. By the quotient property of the
projective tensor product, the operator
$q_\mu\widehat{\otimes}_\pi q_\nu:
\F(G)\widehat{\otimes}_\pi\F(H)
\longrightarrow
\F_\mu^h(G)\widehat{\otimes}_\pi\F_\nu^h(H)$
is a quotient map whose kernel is
$K_{\mu,\nu}$.
Hence
\[
\bigslant{
\F(G)\widehat{\otimes}_\pi\F(H)
}{
K_{\mu,\nu}
}
\cong
\F_\mu^h(G)\widehat{\otimes}_\pi\F_\nu^h(H)
\]
canonically and isometrically. Combining this with
\eqref{eq:HBLip-quotient-predual} gives
\eqref{eq:harmonic-tensor-factorization}.
\end{proof}


\subsubsection{Finite-dimensional consequences and further questions}

The preceding tensor-predual representation yields an immediate
finite-dimensional consequence. Recall that
$\LHFz^{\mu}(G)
=
\{f\in\Lip_0(G):P_\mu f=f\}$
and, by Proposition~\ref{prop:one-variable-harmonic-predual},
$\LHFz^{\mu}(G)
\cong
\bigl(\F_\mu^h(G)\bigr)^*$. The analogous identification holds for $(H,\nu)$.

\begin{corollary}\label{cor:dimension-harmonic-two-lipschitz}
Let $G$ and $H$ be finitely generated groups and let $\mu$ and $\nu$
be finitely supported probability measures on $G$ and $H$, respectively.
Suppose that
$\dim \LHFz^{\mu}(G)=m<\infty
~\text{and}$ $
\dim \LHFz^{\nu}(H)=n<\infty$ (such finite dimensional space exists, see \cite{DouchaIPC}).
Then
$\dim \HBLipz^{\mu,\nu}(G,H)=mn$.
\end{corollary}

\begin{proof}
By Proposition~\ref{prop:one-variable-harmonic-predual},
$\LHFz^{\mu}(G)
\cong
\bigl(\F_\mu^h(G)\bigr)^*$
and
$\LHFz^{\nu}(H)
\cong
\bigl(\F_\nu^h(H)\bigr)^*$. Since the spaces on the left are finite-dimensional, we
obtain
$\dim \F_\mu^h(G)=m,
~\dim \F_\nu^h(H)=n.$
For finite-dimensional Banach spaces, the completed projective tensor
product coincides algebraically with the ordinary tensor product.
Consequently,
we have $\dim\bigl(
\F_\mu^h(G)\widehat\otimes_\pi\F_\nu^h(H)
\bigr)
=
mn$.
 Since a finite-dimensional vector space and its dual
have the same dimension, it follows that
$\dim \HBLipz^{\mu,\nu}(G,H)=mn$ by the Proposition~\ref{prop:harmonic-tensor-predual},
$\HBLipz^{\mu,\nu}(G,H)
\cong
\bigl(
\F_\mu^h(G)\widehat\otimes_\pi\F_\nu^h(H)
\bigr)^*$.
\end{proof}

If $\dim\LHFz^\mu(G)<\infty$, then
$\F_\mu^h(G)$ is finite-dimensional, and hence
$N_\mu(G)$ is complemented in $\F(G)$.
This leads to the following broader question.

\begin{question}\label{ques:Nmu-complemented}
Let $G$ be a finitely generated group and let $\mu$ be a finitely
supported probability measure on $G$. Under what conditions on $G$ and
$\mu$ is
$N_\mu(G)$
complemented in $\F(G)$?
Equivalently, under what conditions does the canonical quotient map
$q_\mu:\F(G)\longrightarrow\F_\mu^h(G)$
admit a bounded linear right inverse?
\end{question}
The next question is related to the above one in two-variable set up.
\begin{question}\label{ques:tensor-quotient-splitting}
Let $G$ and $H$ be finitely generated groups, and let $\mu$ and $\nu$
be finitely supported probability measures on $G$ and $H$,
respectively. Under what conditions does the canonical quotient map
\[
q_\mu\widehat{\otimes}_\pi q_\nu:
\F(G)\widehat{\otimes}_\pi\F(H)
\longrightarrow
\F_\mu^h(G)\widehat{\otimes}_\pi\F_\nu^h(H)
\]
admit a bounded linear right inverse?
\end{question}

A positive answer to Question~\ref{ques:tensor-quotient-splitting}
would yield a bounded projection onto the corresponding harmonic
two-Lipschitz space. Indeed, if
$S:
\F_\mu^h(G)\widehat{\otimes}_\pi\F_\nu^h(H)
\longrightarrow
\F(G)\widehat{\otimes}_\pi\F(H)$
is a bounded right inverse of
$q_\mu\widehat{\otimes}_\pi q_\nu$, then for
$T\in\BLip_0(G,H)$ one may define $\mathcal PT$ through
$(\mathcal PT)_L
=
(T_L\circ S)
\circ
(q_\mu\widehat{\otimes}_\pi q_\nu)$.
This gives a bounded projection
\[
\mathcal P:
\BLip_0(G,H)
\longrightarrow
\HBLipz^{\mu,\nu}(G,H),
\qquad
\|\mathcal P\|\leq \|S\|.
\]
The two questions are related. If $q_\mu$ and $q_\nu$ admit bounded
linear right inverses $S_\mu$ and $S_\nu$, respectively, then
$S_\mu\widehat{\otimes}_\pi S_\nu$
is a bounded linear right inverse of
$q_\mu\widehat{\otimes}_\pi q_\nu$.
Thus complementability of the one-variable harmonic relation spaces
gives a sufficient condition for complementability in the
two-variable setting.

\begin{question}\label{ques:separate-product-harmonicity}
Let $\mu$ and $\nu$ be finitely supported probability measures on
finitely generated groups $G$ and $H$, respectively. Under what
conditions does harmonicity with respect to the product measure
$\mu\otimes\nu$ imply harmonicity in each variable? Equivalently, when
does
\[
\ker(P_{\mu\otimes\nu}-I)\cap\BLip_0(G,H)
=
\ker(P_\mu^{(1)}-I)
\cap
\ker(P_\nu^{(2)}-I)
\cap
\BLip_0(G,H).
\]
\end{question}

\noindent\textbf{Declaration of generative AI and AI-assisted technologies in the manuscript preparation process.}
The ideas developed in this paper originated solely with the author. It was initiated during his doctoral studies. During the final stage of manuscript preparation, the author used OpenAI's ChatGPT, accessed through ChatGPT Plus, to improve the English, obtain suggestions for references, clarify standard facts about word metrics on finitely generated groups, and organise the presentation of the material. 
Some of the results also obtained by repetative chatting with AI by the author.
The author independently verified the mathematical statements and proofs, reviewed and revised the AI-assisted content, and takes full responsibility for the final version of the article.

\medskip
\noindent\textbf{Acknowledgements.}
This work was initiated while the author was a research scholar at the
School of Mathematical Sciences, National Institute of Science Education
and Research (NISER), Homi Bhabha National Institute, India.
The author thanks Professor Michal Doucha for suggesting relevant references during the development of this work. Some of the author's ideas were also shared with him through personal correspondence. The author also thanks Professor Pablo Turco for discussions during his doctoral studies concerning certain results on Lipschitz spaces. These results are adapted here to the two-variable setting, particularly in Remark~\ref{rem:general-codomain-bidual}.

\medskip
\noindent The author declares that he has no competing interests.

\medskip
\noindent\textbf{Data availability.}
No datasets were generated or analysed during this study.


\begin{thebibliography}{99}

\bibitem{LHF}
G.~Amir et al.,
Lipschitz harmonic functions on vertex-transitive graphs,
\emph{Electron. Commun. Probab.} \textbf{29} (2024),
Paper No.~21, 4 pp.; MR4756365.

\bibitem{OEUATS}
R.~F. Arens and J.~Eells,
On embedding uniform and topological spaces,
\emph{Pacific J. Math.} \textbf{6} (1956), 397--403;
MR0081458.

\bibitem{GNFA}
Y.~Benyamini and J.~Lindenstrauss,
\emph{Geometric nonlinear functional analysis. Vol. 1},
American Mathematical Society Colloquium Publications, 48,
American Mathematical Society, Providence, RI, 2000;
MR1727673.

\bibitem{PILFSIBGA}
M.~C\'uth and M.~Doucha,
Projections in Lipschitz-free spaces induced by group actions,
\emph{Math. Nachr.} \textbf{296} (2023), no.~8, 3301--3317;
MR4626883.

\bibitem{OTSOLFS}
M.~C\'uth, M.~Doucha and P.~Wojtaszczyk,
On the structure of Lipschitz-free spaces,
\emph{Proc. Amer. Math. Soc.} \textbf{144} (2016), no.~9,
3833--3846; MR3513542.

\bibitem{DouchaIPC}
M.~Doucha,
\emph{Invariant pointwise closed subspaces of Lipschitz spaces and their preduals},
arXiv:2608.24143 [math.FA], 2026.

\bibitem{DAS}
M.~J. Field,
\emph{Dynamics and symmetry},
ICP Advanced Texts in Mathematics, 3,
Imp. Coll. Press, London, 2007; MR2358722.

\bibitem{LFBS}
G.~Godefroy and N.~J. Kalton,
Lipschitz-free Banach spaces,
\emph{Studia Math.} \textbf{159} (2003), no.~1, 121--141;
MR2030906.

\bibitem{MSLOGFHM}
F.~Gu\'eritaud and F.~Kassel,
Maximally stretched laminations on geometrically finite hyperbolic manifolds,
\emph{Geom. Topol.} \textbf{21} (2017), no.~2, 693--840;
MR3626591.

\bibitem{TLOI}
K.~Hamidi et al.,
Two-Lipschitz operator ideals,
\emph{J. Math. Anal. Appl.} \textbf{491} (2020), no.~2,
124346, 21 pp.; MR4122068.

\bibitem{AM3}
A.~K. Karn and A.~Mandal,
\emph{Position of $L(X,Y)$ in $\operatorname{Lip}_0(X,Y)$},
preprint, arXiv:2506.09324 [math.FA], 2025.

\bibitem{EK} E. Kowalski, {\it An introduction to the representation theory of groups}, Graduate Studies in Mathematics, 155, Amer. Math. Soc., Providence, RI, 2014; MR3236265

\bibitem{AM2}
A.~Mandal,
A proper closed subspace of the Lipschitz dual containing the linear dual,
\emph{Monatsh. Math.} (2026),
\url{https://doi.org/10.1007/s00605-026-02223-0}.

\bibitem{AMA}
A.~Mandal,
\emph{Group Actions, Fixed Points and Orbit Quotients in Lipschitz Spaces},
arXiv:2609.25931 [math.FA], 2026.


\bibitem{AM1}
A.~Mandal,
Projection from space of two-Lipschitz operators onto the space of
bilinear maps,
\emph{Mediterr. J. Math.} \textbf{23} (2026), no.~2,
Paper No.~75, 17 pp.; MR5038008.

\bibitem{LFSATAGA}
M.~G. Megrelishvili,
Lipschitz-free spaces: a topometric approach and group actions,
\emph{Mediterr. J. Math.} \textbf{22} (2025), no.~8,
Paper No.~202, 39 pp.; MR4973079.


\bibitem{RAR}
R.~A. Ryan,
\emph{Introduction to tensor products of Banach spaces},
Springer Monographs in Mathematics,
Springer, London, 2002; MR1888309.

\bibitem{SerreRep}
J.-P.~Serre,
\emph{Linear Representations of Finite Groups},
Graduate Texts in Mathematics, Vol.~42,
Springer-Verlag, New York, 1977.

\bibitem{JLT}
J.~L. Taylor,
A joint spectrum for several commuting operators,
\emph{J. Functional Analysis} \textbf{6} (1970), 172--191;
MR0268706.

\bibitem{tomDieck}
T.~tom Dieck,
\emph{Representation Theory},
Preliminary version, 9 February 2009,
Georg-August-Universit\"at G\"ottingen,
available at
\url{https://studyres.com/doc/6233225/representation-theory-tammo-tom-dieck}.

\bibitem{LA}
N.~Weaver,
\emph{Lipschitz algebras}, second edition,
World Scientific, Hackensack, NJ, 2018;
MR3792558.

\bibitem{Yadin}
A.~Yadin,
\emph{Harmonic functions and random walks on groups},
Cambridge Studies in Advanced Mathematics, 213,
Cambridge Univ. Press, Cambridge, 2024;
MR4742316.

\end{thebibliography}
\end{document}